\documentclass[11pt]{amsart}
\usepackage{amscd,amsfonts,amsmath,amssymb,amstext,amsthm, xcolor}
\usepackage{makecell}
\usepackage[pdftex]{graphicx}
\usepackage{graphicx}
\usepackage{caption}
\usepackage[linktocpage=true]{hyperref}
\usepackage{lscape}
\usepackage{multirow}
\usepackage{tikz}
\usepackage{tikz-cd}
\usetikzlibrary{matrix,arrows}
\usepackage{geometry}
\usetikzlibrary{matrix,positioning,decorations.pathreplacing}
\usepackage{ytableau}
\date{\today}
\theoremstyle{plain}
\newtheorem{theorem} {Theorem}
[section]

\newtheorem{proposition}[theorem]{Proposition}
\newtheorem{corollary} [theorem]{Corollary}
\theoremstyle{definition}

\newtheorem*{example} {Example}
\numberwithin{equation}{section}
 \def\ind{\operatorname{ind}\,}
  \def\ext{\operatorname{ext}}
\newcommand{\ve}{\varepsilon}

\title[Ampleness of normal bundles of base cycles in flag domains]{Ampleness of normal bundles of base cycles in flag domains}

\author[J. Hong]{Jaehyun Hong}
 \address{Center for Complex Geometry,  Institute for Basic Science (IBS), Daejeon 34126, Republic of Korea}
    \email{jhhong00@ibs.re.kr}

\author[A. Seo]{Aeryeong Seo}
\address{Department of Mathematics, Kyungpook National University, Daegu, 41566, Republic of Korea}
\email{aeryeong.seo@knu.ac.kr}


\begin{document}

\begin{abstract} 
%Flag manifolds are closed orbits of complex semisimple Lie groups. Their noncompact real forms have only finitely many open orbits, called flag domains. 
Flag domains    are   open  orbits of noncompact real forms   of complex semisimple Lie groups    acting on  flag manifolds. To each flag domain  one can associate a compact complex manifold called the base cycle.
The ampleness of the normal bundle of the base cycle in a flag domain    measures the concavity near the base cycle.
In this paper we compute the ampleness of normal bundles of base cycles in flag domains in various cases, including flag domains in the full flag manifolds $G/B$ when $G$ is classical, and period domains  parameterizing polarized Hodge structures with   fixed Hodge numbers.
\end{abstract}
\maketitle

\def\Label#1{\label{#1}{\bf (#1)}~}

\section{Introduction}

By flag domain we mean an open  orbit of a noncompact real form $G_0$ of a complex semisimple Lie group $G$ acting on a flag manifold $G/P$.
%The most two well known flag domains are Hermitian symmetric spaces of the noncompact type and period domains.
To each flag domain $D$ it is possible to associate a compact homogeneous complex manifold called the base cycle in $D$. It is $K$-orbit in $D$ where $K$ is the complexification of a maximal compact subgroup in $G_0$. One of the most studied flag domains are Hermitian symmetric spaces of noncompact type. These domains are pseudoconvex, in a sense that there exists a
plurisubharmonic exhaustion function on $D$, and the base cycles are nothing but points.
On the other hand, most of period domains  parameterizing polarized Hodge structures with   fixed Hodge numbers are flag domains which are not pseudoconvex since their base cycles are of positive dimension. Moreover, these domains are cycle connected, that is, any two points in the domain can be connected by a chain of translates of the base cycle. The geometry of flag domains is deeply related to that of cycle spaces. For the further study, see \cite{FHW}, \cite{H13} and the references therein.

The cycle connectivity of flag domains can be interpreted as the finiteness of  $H^0(D, \mathcal F)$ for any coherent sheaf $\mathcal F$ {as studied} in \cite{HHL} and \cite{HHS}. In \cite{HHS} the authors used  the ampleness formula of Snow \cite{S1, S2} and related  it with the signature of  the Levi curvature of the dual of the normal bundle of the base cycle.

The finiteness of $H^q(D, \mathcal F)$ for any coherent sheaf on a flag domain $D$ is known only for $q=0$ (\cite{HHL}, \cite{HHS}) and it remains open for higher $q$. According to the Andreotti-Grauert finiteness theorem, the finiteness of $H^q(D, \mathcal F)$ is determined by the degree of concavity of $D$: $H^i(D, \mathcal F)$ is finite for $i \leq  dih(\mathcal F)-q$ if $q>0$ and $D$ is $q$-concave, i.e., there is an exhaustion function on $D$ which is $q$-concave outside a compact set, where $dih(\mathcal F)$ is the homological dimension of $\mathcal F$.  However the degrees of concavity of flag domains are not known.

In this paper, we consider a weaker notion of concavity. Given a compact complex submanifold $X$ of a complex manifold $Y$, we say that $Y$ is $q$-concave near $X$ if there are an open neighborhood $U$ of $X$ in $Y$ and an exhaustion function $\rho$ on $U$ which is $q$-concave outside $X$.   % and a smooth function $\rho:U \rightarrow \mathbb R^{\geq 0}$ with $\{\rho=0\}=X$ and $d \rho \not=0$ otherwise, so that for $\ve >0$ sufficiently small $U_{\ve}:=\{\rho <\ve\}$ is a relatively compact neighborhood of $X$ and that $\rho $ is $\ell$-concave on $U_{\ve} - X$. In other words, the Levi form of $\rho$ has at most $(\ell-1)$-non-negative eigenvalues.
This weaker concavity of a flag domain $D$ is related with the ampleness of the normal bundle of the base cycle $C$ of $D$ as follows.
 A flag domain $D$ is $(r+a +1)$-concave near $C$, where $r$ is the codimension of $C$ in $D$ and $a:=a(N_{C|D})$ is the ampleness of the normal bundle of $C$ in $D$ (Corollary 3 of \cite{HHS}). In other words, there is an exhaustion function $\rho$ near $C$ such that the Levi form has at least $\dim_{\mathbb C} C - a $ negative eigenvalues.  %(Corollary 4 of \cite{HHS}).
{Furthermore, the pseudoconvexity of $D$ can be described as follows.
 \begin{proposition} [Proposition 10 \cite{HHS}] \label{pseudoconvex}
     If $C$ is a closed $K$-orbit in $G/P$ in a flag domain $D$, then $a= \dim_{\mathbb C}C$ if and only if the normal bundle $N_{C|D}$ is trivial and $D$ fibers over a Hermitian symmetric domain.
 \end{proposition}}
In this regard, it is worth computing the ampleness $a(N_{C|D})$ of the normal bundle of $C$ in $D$, or the difference $\dim_{\mathbb C}C - a(N_{C|D})$,   which measures the concavity near the base cycle when $D$ is not pseudoconvex.

\medskip

Before stating our theorem, we introduce some notations. 
Let $p$, $q$ be positive integers. Consider the set $\mathcal P_{p,q}$ of all subsets of cardinality $p$ in $\{1,2, \dots, p+q\}$. 
Define functions
 $$I^+: \mathcal P_{p,q} \rightarrow \mathbb N \text{ and } I^-:\mathcal P_{p,q} \rightarrow \mathbb N$$    as follows. 
For a subset ${\bf j}=\{j_1, \dots, j_p\}$ of $\{1, \dots, p+q\}$
with $j_1 <j_2 < \dots <j_p$,
let $\{j_{p+1},\ldots, j_{p+q}\}$ denote the complement $\{1,2, \dots,{p+q} \} \backslash {\bf j}$ with 
$j_{p+1}<j_{p+2}<\cdots<j_{p+q}$.

Define $h^+({\bf j})$ and $h^-({\bf j})$ by 
%Let $(a', b') $ and $(a'', b'')$ with $1\leq a', a''\leq p<p+1\leq b', b''\leq p+q$ be such that
\begin{equation}\label{ab_int}
\begin{aligned}
 h^+({\bf j})&:=\left\{
\begin{array}{cc}
\min\{ b-a: a \leq p <p+1 \leq b, j_a <j_b \}&\text{ if }\, {\bf j} \neq \{q+1, \ldots, p+q\},\\
0&\text{ if } \,{\bf j} = \{q+1, \ldots, p+q\},\\
\end{array}\right.\\
h^-({\bf j})&:=\left\{
\begin{array}{cc}
     {\max}\{ b-a: a \leq p <p+1 \leq b, j_a >j_b\} 
     &\text{ if }\, {\bf j} \neq \{1,\ldots, p\},\\
     0 
     &\text{ if }\, {\bf j} = \{1,\ldots, p\},\\
\end{array}\right.
\end{aligned}
\end{equation}
respectively. 
%\begin{equation}\label{ab_int}
%\begin{aligned}
% b'-a'
% h^+({\bf j})&:=\min\{ b-a: 1\leq a\leq p<   b\leq p+q, \,\, j_a <j_b \}, \text{ and }\\
% a'' - b''
% h^-({\bf j})&:=\max\{ b-a:1\leq a\leq p<   b\leq p+q,\,\, j_a >j_b\}
%\end{aligned}
%\end{equation}
\noindent For a combinatorial meaning of 
the quantities $h^+({\bf j}) $ and $h^-({\bf j}) $, see the end of this section. 
Define $I^+({\bf j})$ and $I^-({\bf j})$ by
\begin{equation}\nonumber
\begin{aligned}
%I^+({\bf j}) &:= (b'-(p+1)) + (p-a') +p, \text{ and }\\
%I^-({\bf j}) &:= (a''-1)+ (p+q-b'') +q 
I^+({\bf j}) &:=h^+({\bf j}) -1 +p, \text{ and }\\
I^-({\bf j}) &:=(p+q) - h^-({\bf j})-1 +q
\end{aligned}
\end{equation}
respectively. 

If the minimum occurs at $(a',b')$ (the maximum occurs at $(a'', b'')$, respectively) in the right hand side of  \eqref{ab_int}, then 
we get 

\begin{equation}\nonumber
\begin{aligned}
 I^+({\bf j}) & = (b'-(p+1)) + (p-a') +p, \text{ and }\\
 I^-({\bf j}) & = (a''-1)+ (p+q-b'') +q.  
%I^+({\bf j}) &:=h_{\min}({\bf j}) -1 +p, \text{ and }\\
%I^-({\bf j}) &:=(p+q) - h_{\max}({\bf j})-1 +q
\end{aligned}
\end{equation}
%respectively. 

\medskip 
Our first result is the following.

\begin{theorem} [Section \ref{sect:full flag manifolds}] Let $G_0$ be a classical simple Lie group with its Lie algebra $\frak g_0$ and its complexification $G$. Let $G/B$ be a
flag manifold where $B$ is a Borel subgroup of $G$, and let $D$ be a flag domain in $G/B$. Let $a(N_{C|D})$ denote the ampleness of the normal bundle of the base cycle $C$ in  $D$. 
%Then for each base cycle $C$ in $G/B$, the ampleness of its normal bundle is $\dim_{\mathbb C} C - \ind (-\Lambda_{\text{ext}}(E_0))$ where $\ind (-\Lambda_{\text{ext}}(E_0))$ 
Then $\dim_{\mathbb C} C -a(N_{C|D})$ 
is given as follows:\\

\begin{center}
\begin{tabular}{c|c|c}
\hline \\[-8pt]
$\frak g_0$
& $W^\theta_1$
& $ \dim_{\mathbb C} C -a(N_{C|D})$ \\[5 pt]
%& $\ind (-\Lambda_{\text{ext}}(E_0))$ \\[5 pt]
\hline && \\[-8pt]
$\frak{sl}(m,\mathbb R)$
& {\makecell[c] {
$\{\text{id}\, \}$ if $m=2r+1$, $r\geq 1$\\
$\{\text{id}, s_{\gamma_r}\}$ if $m=2r$, $r\geq 1$}}
&$\left[ m-1 \over 2 \right]$\\[5 pt]
\hline && \\[-8pt]
$\frak{su}(p,q)$, $p\geq 2$
&$\{w_{\bf j} : {\bf j}\subset \{1,\ldots, p+q\}, |{\bf j}|=p\}$
&{$\min\{ I^+({\bf j}) -p,I^-({\bf j}) -q\}$}
\\[5 pt]
\hline && \\[-8pt]
$\frak{su}(1,q)$
&$\{w_{j} : {j}\in\{0,1,\ldots, q\}\}$
&$\min\{ j, q-j-1\}$\\[5 pt]
\hline && \\[-8pt]
$\frak{sp}(r,\mathbb R)$
&$\{w_{\bf j} : {\bf j}\subset \{1,\ldots, r\}\}$
&$\min\{ |{\bf j}|, r-|{\bf j}|\}$\\[5 pt]
\hline && \\[-8pt]
\makecell[c]{$\frak{so}(2p+1, 2q+1)$, $1\leq p$\\ or $\frak{sp}(p,q)$  $p\leq q$}
&$\{w_{\bf j} : {\bf j}\subset \{1,\ldots, p+q\}, |{\bf j}|=p\}$
&$\min \{ I^+({\bf j}), I^-({\bf j}) \}$\\[5 pt]
\hline && \\[-8pt]
{$\frak{so}(1,2q+1)$}
&{$\{id\}$}
&{$q$}\\[5 pt]
\hline && \\[-8pt]
\makecell[c]{$\frak{so}(2p, 2q+1)$, \\$p\geq 3$, $ q\geq 1$}
&$\{w_{\bf j}, w_{\bf j}' : {\bf j}\subset \{1,\ldots, p+q\}, |{\bf j}|=p\}$
&$\min \{ I^+({\bf j})-1, I^-({\bf j}) \}$\\[5 pt]
\hline && \\[-8pt]
\makecell[c]{$\frak{so}(2p, 2q)$,\\ $2\leq p\leq q$}
&$\{w_{\bf j}, w_{\bf j}' : {\bf j}\subset \{1,\ldots, p+q\}, |{\bf j}|=p\}$
&$\min \{ I^+({\bf j})-1, I^-({\bf j}) -1\}$\\[5 pt]
\hline && \\[-8pt]
$\frak{so}(2, q)$
&$\{w_{j} : {j}\in\{1,\ldots, q+1\}\}$
&$j-1$\\[5 pt]
\hline && \\[-8pt]
$\frak{so}(4, q)$
&$\{w_{jk} : 1\leq j<k\leq q+1\}$
&
\makecell[c]{$j+1$  if $j+1\neq k$,\\
$j\quad \,$ if $j+1=k$} \\
[5 pt]
\hline && \\[-8pt]
$\frak{sl}(m, \mathbb H)$
& {$\{ \text{id}\,\}$}
&$m$\\[5 pt]
\hline
\end{tabular}

\end{center}
Here, the set of base cycles is parameterized by $W_1^{\theta}$ (For notations, see Section \ref{sect:classification of base cycles} {and Section~\ref{sect:full flag manifolds}}).

\end{theorem}
The same method allows one to compute the ampleness associated to every base cycle of the quotient $G/Q$, with $Q$ a parabolic subgroup of $G$. {In what follows} we will restrict our attention to the base cycles of period domains, which are important examples of flag domains.
%For a given general parabolic subgroup $Q$ in $G$, for base cycles in $G/Q$.
%Instead of computing ampleness of every base cycle, we do this for the base cycles in period domains, important examples of flag domains.

\medskip

Let $D$ be a period domain parameterizing polarized Hodge structures with fixed Hodge numbers $\{h^{r,s}\}_{r+s=n}$.  %For $0 \leq p \leq n$, set $f^p:=\sum_{s=p}^n h^{s, n-s} = \sum_{s=0}^{n-p} h^{n-s,s}$. Set $f:=f^0$.
{Set 
$ f^r:=   h^{n,0} + h^{n-1,1}+ \cdots   + h^{r,n-r} $ for $0 \leq r \leq n$.  } 

\medskip
If $n=2k+1$ for some $k\in \mathbb N$, then
$$D=Sp(m,\mathbb R)/U(h^{n,0})\times U(h^{ n-1,1}) \times \dots \times U(h^{ k+1,k}), $$
where $ m={\displaystyle \sum_{r\geq k+1}h^{r,s}}$.
Let
\begin{equation}\nonumber
h_o:={\displaystyle \sum_{\substack{ r\geq k+1 \\ r, odd}}h^{r,s}} \quad\text{ and }\quad
h_e:={\displaystyle \sum_{\substack{r \geq k+1 \\ r, even}}h^{r,s}}.
\end{equation}

\medskip
If $n=2k$ for some $k\in \mathbb N$, then
$$D=SO(m_e, m_o)/U(h^{n,0})\times U(h^{ n-1,1}) \times \dots \times U(h^{ k+1, k-1 })\times SO(h^{k,k}), $$
where $m_e:=\sum_{r, even} h^{r,s}$ and $m_o:=\sum_{r, odd}h^{r,s}$. %{ Set $m:=\left[\frac{  m_e + m_o }{ 2}\right]$. } 
% $m:=\left[ m_e + m_o  \over 2\right]$.

Let 
\begin{equation}\nonumber
\begin{aligned}
&p=\frac{m_e}{2},\,
\ell =m_o, \, q=\left[m_o \over 2\right] 
\text{ and }\\&
{\bf j} = {\{f^{n}+1, \dots, f^{n-1}, \quad \dots, \quad f^{k+1}+1, \ldots, p+q}\} \subset \{1, \dots, { p+q}\} \quad \text{ if } k \text{  is odd},
\end{aligned}
\end{equation} 
 and let
\begin{equation}\nonumber
\begin{aligned}
&p=\frac{m_o}{2},\,
\ell =m_e, \, q=\left[m_e \over 2\right] 
\text{ and } \\
&{\bf j} = { \{1, \dots, f^{n},  \quad \dots, \quad    f^{k+1}+1, \ldots, p+q}\}  \subset \{1, \dots, { p+q}\} \quad  \text{ if } k \text{  is even}. 
\end{aligned}
\end{equation}

%For each $k$ let $(a', b')$, $(a'', b'')$ be as in
%\eqref{ab_int}.
Our second result is the following.
\begin{theorem} [Section \ref{sect:general flag manifolds}]
Let $D$ be a period domain parameterizing polarized Hodge structures with fixed Hodge numbers $\{h^{r,s}\}_{r+s=n}$.
Let $a(N_{C|D})$ denote the ampleness of the normal bundle of the base cycle $C$ in $D$. Then 
$\dim_{\mathbb C}C -a(N_{C|D})$  is given as follows:\\
%The ampleness of the normal bundle of the base cycle $C$ in $D$ is %$\dim_{\mathbb C} - \ind(-\Lambda_{\text{ext}}(E_0))$ where
%$\ind(-\Lambda_{\text{ext}}(E_0))$ is given as follows:\\
\begin{center}
\begin{tabular}{c|c|c}
\hline && \\[-8pt]
weight $n$
&$\frak g_0$
& $\dim_{\mathbb C} C -a(N_{C|D})$ \\[5 pt] 
%& $\ind (-\Lambda_{\text{ext}}(E_0))$ \\[5 pt]
\hline && \\[-8pt]
{$n=2k+1$}
&{$\frak{sp}(m, \mathbb R)$}
&$\min\{ h_o, h_e\}$\\[5pt]
\hline && \\[-8pt]
{$n=2k$}
&{$\frak{so}(2p, \ell)$}
& \makecell[c] { $\min\{I^+({\bf j})-1, \,I^-({\bf j}) \}$ if $\ell$ is odd \\
$\min\{I^+({\bf j})-1, \,I^-({\bf j}) -1\}$ if $\ell$ is even} 
%$\min\{I^+(a',b')-1, \,I^-(a'',b'')-\delta_{\ell,e}\}$
\\[10 pt]
\hline
\end{tabular}

\end{center}
%where
%and $\delta_{\ell, e}$ is $1$ if $\ell$ is even and $0$ otherwise.
\end{theorem}

\medskip 
{
 It is worth mentioning a combinatorial meaning of 
{the} quantities $h^+({\bf j})=b'-a'$ and $h^-({\bf j})=b''-a''$ introduced in \eqref{ab_int}.  
Consider the Young diagram of shape $p^q=\underbrace{(p, \dots, p)}_{q}$. The hook length at a box $u$ is defined by the number of boxes directly to the right or directly below $u$, counting $u$ itself once. 
We label columns  with  $j_1, \dots, j_p$ from the left to the right, and rows with $j_{p+1}, \dots, j_{p+q}$ from the bottom to the top.   
This labeling defines   coordinates $(j_a,j_b)$ on boxes. % in such a way that the first coordinate is given by $j_1, \dots, j_p$ from the left to the right, and the second coordinate is given by $j_{p+1}, \dots, j_{p+q}$ from the bottom to the top.   
Then the hook length at a box with coordinate $(j_a,j_b)$ is $b-a$. 
%First consider $p\times q$ boxes, which are labelled with $j_1,\ldots, j_p$ on the horizontal side and with $j_{p+1},\ldots, j_{p+q}$ on the vertical side. 
%Let us remind the reader that the hook length of a box whose position is $(j_a, j_b)$ is defined as the number of boxes of the position $(j_\alpha, j_\beta)$ such that $\alpha\geq a$ and $\beta\leq b$.
%Next

For example, consider the set $\{1,2,\ldots, 7\}$ and ${\bf j} = \{2,5,6\}$.
In this case one has $j_1 = 2$, $j_2=5$, $j_3=6$ and $j_4=1$, $j_5 = 3$, $j_6=4$, $j_7=7$.
Then the $3\times 4$ boxes are given as in the following diagram.
The numbers written in the boxes are the corresponding hook lengths.

\begin{center}
\ytableausetup{centertableaux}
\begin{ytableau}
\none[j_7] &  6 &  5 &  4 \\
\none[j_6] &   5& 4 & 3 \\
\none[j_5] &  4& 3 & 2 \\
\none[j_4]  & 3 &2  & 1 \\
\none & \none[j_1] & \none[j_2] & \none[j_3]
\end{ytableau}
\end{center}

Color those boxes whose coordinate $(j_a, j_b)$ satisfies $j_a<j_b$. Then the minimum hook length of   colored boxes is equal to ${h^+({\bf j})=b'-a'}$, while the maximum hook length of   uncolored boxes is given by ${h^-({\bf j})=b''-a''}$.

For example, from the coloring of the boxes with $j_a<j_b$ 
\begin{center}
\ytableausetup{centertableaux}
\begin{ytableau}
\none[7] & *(gray!40)6 & *(gray!40)5 & *(gray!40)4 \\
\none[4] & *(gray!40) 5& 4 & 3 \\
\none[3] & *(gray!40)4& 3 & 2 \\
\none[1]  & 3 &2  & 1 \\
\none & \none[2] & \none[5] & \none[6]
\end{ytableau}
\end{center}
we get $b'-a' = 4$ and $b''-a''= 4$.

}

\medskip

For the computation, we use Snow's ampleness formula for homogeneous vector bundles $E=L\times_PE_0$ on a flag manifold $C=L/P$   (Proposition~\ref{Snow ampleness formula}):
$$a(E) = \dim_{\mathbb C} C -  \ind  (-\Lambda_{\ext}(E_0)).$$
Here,  the index $\ind  (-\Lambda_{\ext}(E_0))$ is an invariant of $E_0$ as a $P$-module (For the definition, see Section \ref{sect: ampleness formula}).

\medskip

The paper is organized as follows. In Section \ref{sect: ampleness} we recall  the definition and properties of ampleness.  Our computation is based on the results in Chapter 17--18 of \cite{FHW}.  We review their  theory on classifications of base cycles in Section \ref{sect:classification of base cycles} and describe extremal weights in $E_0$ in  Section  \ref{sect:extremal weights and indices}. We summarize our method how to compute the ampleness of the normal bundle of the base cycle  in Section \ref{sect:our strategy}.
Following this strategy we compute the ampleness of the normal bundles of the base cycles of flag domains  in Section \ref{sect:full flag manifolds}   and  we do this for period domains in Section \ref{sect:general flag manifolds}.

{\bf Acknowledgement}
The authors would like to express their gratitude to the anonymous referee for carefully reading the manuscript and providing excellent suggestions for improvement.
The first author was supported by the Institute for Basic Science (IBS-R032-D1). The second author was partially supported by Basic Science Research Program through the National Research Foundation of Korea (NRF) funded by the Ministry of Education (NRF-2022R1F1A1063038).

\section{Ampleness} \label{sect: ampleness}

\subsection{Definition}

%Let $E$ be a holomorphic vector bundle on a complex manifold $X$.
%We say that $E$ is {\it $k$-ample} if some power $\zeta^m$ of the hyperplane line bundle $\zeta$ on $\mathbb P(E^*)$ is spanned and   the induced map $$\mathbb P\nu_m:\mathbb P(E^*) \rightarrow \mathbb P(H^0(\mathbb P(E^*), \zeta^m)^*)$$ has fiber dimension at most $k$.
   %By Sommese \cite{So78}, the maximum fiber dimension does not depend on a choice of $m$ among all $m$ for which $\zeta^m$ is spanned.

  Let $L$ be a line bundle on a complex manifold $X$. We say that $L$ is

  \begin{enumerate}
  \item {\it $k$-ample} if some power of $L$ is spanned and the induced map
  $$\mathbb P\nu_m:M \rightarrow \mathbb P(H^0(M, L^m)^*)$$
  has fiber dimension at most $k$;

  \item  {\it $k$-positive} if there is a Hermitian metric on $L$ whose curvature has at least $(n-k)$ positive eigenvalues;

  \item  {\it cohomologically $k$-ample} if for every coherent sheaf $\mathcal F$  on $M$, there is $m_0>0$ depending on $\mathcal F$ such that
      $$H^i(M, \mathcal F \otimes L^m)=0 \text{ for all } m \geq m_0$$
    for all $i > k$.
    \end{enumerate}

%if some multiple $L^{\otimes m}$ of $L$ is globally generated.

    \begin{proposition} [\cite{S1}, \cite{S2}, \cite{M}] \label{prop:equivalent} Let $L$ be a line bundle on a compact complex manifold such that some power of $L$ is spanned. Then the followings are equivalent.

    \begin{enumerate}
    \item $L$ is $k$-ample
    \item $L$ is $k$-positive
    \item $L$ is cohomologically $k$-ample.
    \end{enumerate}

    \end{proposition}

    Let $E$ be a holomorphic vector bundle on a complex manifold $X$. We say that $E$ is {\it $k$-ample} ({\it $k$-positive}, {\it cohomologically $k$-ample}, respectively) if the hyperplane line bundle $\zeta$ on $\mathbb P(E^*)$ is. If some powers of $\zeta$ is spanned, then all three conditions are equivalent (Proposition~\ref{prop:equivalent}).
    The {\it ampleness} $a(E)$ of $E$  is the minimum $k$ such that $E$ is $k$-ample, i.e., the ampleness  $a(E)$ of $E$ is the maximum fiber dimension of the map $\mathbb P \nu_m:\mathbb P(E^*) \rightarrow \mathbb P(H^0(\mathbb P(E^*), \zeta^m)^*)$.

\subsection{Concave neighborhoods of compact submanifolds}

Let $X$ be a compact complex submanifold of a complex manifold $Y$.
We say that $Y$ is $(r+q+1)$-concave near $X$ if  there is a smooth function $\rho :U\to \mathbb{R}^{\ge 0}$ with minimum $\{\rho =0\}=X$ and $d\rho \not=0$ otherwise
so that for $\ve >0$ sufficiently small $U_\ve =\{\rho<\ve \}$ is a relatively compact neighborhood of $X$
and the Levi form of $\rho$ has at most $(r+q)$ non-negative eigenvalues on $U_\ve \setminus X$.

\begin{proposition} [Corollary 3 of \cite{HHS}] Let $X$ be a compact complex submanifold of a complex manifold $Y$ and let $E$ be its normal bundle. Let $r$ be the rank of $E$. Assume that $E$ is spanned. Then $Y$ is $(r+a(E) +1)$-concave near $X$.

\end{proposition}

\begin{corollary} [Corollary 4 of \cite{HHS}] \label{concavity}
If $a(E)< \dim X$, then  there is an exhaustion $\rho$ such that the Levi form at every point of the boundaries of the sublevel sets $\{\rho < \ve \}$
has $ \dim X -a(E)$  negative eigenvalues.
\end{corollary}

%\subsection{ }

%\subsection{From neighborhoods of base cycles to flag domains}

%\vskip 5 cm

% \newpage

\subsection{Ampleness formula} \label{sect: ampleness formula}

Let $E=G\times _P E_0$ be a homogeneous vector bundle on $X=G/P$ where $P$ is a parabolic subgroup of $G$. Then the ampleness $a(E)$ of $E$ can be expressed in terms of the index of a set of weights as follows.

Let $B$ be the Borel subgroup generated by $T$.   Let $W$ be the Weyl group of $G$ and let $\Lambda$ be the set of weights of $G$ with respect to $T$.
 For a $B$-module $E_0$, denote by $\Lambda_{\max}(E_0)$ the set of maximal weights of $E_0$ and by  $\Lambda_{\ext}(E_0)$  the set of extremal weights of $E_0$, which is defined by
 $$\Lambda_{\ext}(E_0) := W(\Lambda_{\max}(E_0)) \cap \Lambda(E_0).$$
The index $\ind(\lambda)$ of a weight $\lambda$   is defined by
$$ \ind(\lambda):= \min\{ \ell(w): w  \in W, w\lambda  \text{ is dominant}\}. $$
%{ For two weights $\lambda, \mu$ in $W.\lambda_I \simeq W^I$, we define the {\it distance} between $\lambda$ and $\mu$ by the minimum of $\ell(w)$ with $w\lambda = \mu$ and $w\in W$.}  
For two  subsets $A$ and $B$ of  $\Lambda$, we define the {\it distance} between $A$ and $B$  by the minimum of $\ell(w)$ with $wA \cap B \not=\emptyset$ and $w\in W$. Then $\ind(\lambda)$ is the distance between $\{\lambda\}$ and the set $\Lambda^+$ of dominant weights. 

For a subset $A$ of $\Lambda$,  the index $\ind(A)$ of $A$  is defined by the minimum of the index $\ind(\lambda)$ of  $\lambda \in A$. { In other words, $\ind(A)$ is the distance between $A$ and the set $\Lambda^+$ of dominant weights.}

\begin{proposition}[Section 2 of \cite{Sn1} and Theorem 8.3 of \cite{Sn2}] \label{Snow ampleness formula}

Let $E=G\times _P E_0$ be a homogeneous vector bundle on $X=G/P$ where $P$ is a parabolic subgroup of $G$. Assume that some power $\zeta^m$ of the hyperplane line bundle $\zeta$ on $\mathbb P(E^*)$ is spanned. Then the ampleness of $E$ is given by
%\begin{eqnarray*}
%a(E) = \max \{ \ell(\omega): \,\omega \in W, \, \Lambda({E_0^*}^U) \cap \omega \Lambda(E_0^*) \not=\emptyset \} + \dim B - \dim P.
%\end{eqnarray*}

$$a(E) = \dim X -  \ind  (-\Lambda_{\ext}(E_0)). $$

\end{proposition}

\section{Normal bundles of base cycles} \label{sect: normal bundles of base cycles}

%\vskip 20 pt

%notations: $X=G/P$ is the base of $E=G\times _P E_0$ here, $X=G/B$ later on.

\subsection{Classifications of base cycles} \label{sect:classification of base cycles}

We will review results in Chapters 17--18 of \cite{FHW}, necessary for our computation of ampleness of normal bundles of base cycles. \\

Let $G_0$ be a noncompact real form of a complex semisimple Lie group $G$ and $K$ be the complexification of a maximal compact subgroup $K_0$ of $G_0$. Let $\theta$ denote the Cartan involution of $G$ associated to $(G,K)$ and let $\tau$ denote the complex conjugation of $\frak g$ over $\frak g_0$. 
The corresponding decomposition into $\pm 1$-eigenvalues of $d\theta$ are $\frak g_0 = \frak k_0 + \frak s_0$ and $\frak g = \frak k + \frak s$ where $\frak k_0$, $\frak k$ are Lie algebras of $K_0$, $K$, respectively.
A closed $K$-orbit $C=K(z)$ in a $G$-flag manifold $Z=G/Q$ is called a base cycle in $Z$.  \\

To a pair $(C,Z)$ of a $G$-flag manifold $Z=G/Q$ and a base cycle $C=K(z)$ in $Z$, we associate  Lie algebra data:
$$\frak b_{\frak k}  \subset \frak b \subset \frak q \subset \frak g,$$
where $\frak g$ is the Lie algebra of $G$, 
$\frak q$ is the Lie algebra of the isotropy group $Q_z$ of $G$ at $z \in Z$, 
$\frak b$ is the Lie algebra of a $\theta$-stable Borel subgroup $B$ of $G$ contained in $Q_z$, and
$\frak b_{\frak k} $ is the Lie algebra of the Borel subgroup $B_{\frak k}  =B \cap K$ of $K$. After we fix a Borel subalgebra $\frak b_{\frak k}^{\rm ref}$ of $\frak k$,  such a Lie algebra data  $\frak b_{\frak k}^{\rm ref} \subset \frak b \subset \frak q \subset \frak g $ determines a base cycle $C$ in a $G$-flag manifold $Z$ uniquely (Section 17.2 of [FHW]). \\

 The space of all {$\theta$-stable} Borel subalgebras $\frak b$ of $\frak g$ containing $\frak b_{\frak k}^{\rm ref}$ can be parameterized  as follows. Let $\frak h$ be a $\tau$-stable $\theta$-stable  Cartan subalgebra of $\frak g$ such that $\frak h_0:=\frak h \cap \frak g_0$ is a maximally compact Cartan subalgebra of $\frak g_0$, i.e., $\frak t:=\frak h_{\frak k}=\frak h \cap \frak k$ is a Cartan subalgebra of $\frak k$.
The Weyl group $W(\frak k, \frak t)$ acts on the centralizer $\frak h$ of $\frak t$ and thus can be considered as a subgroup of the Weyl group $W(\frak g, \frak h)$. Define $W^{\theta}(\frak g, \frak h)$ by the subgroup of $W(\frak g, \frak h)$ consisting of $w \in W(\frak g, \frak h)$ such that $w \circ \theta = \theta \circ w$. Here we consider $w \in W(\frak g, \frak h)$ as an element in $GL(\frak h)$. Fix a Borel subalgebra $\frak b_{\frak k}$ of $\frak k$ and a Borel subalgebra $\frak b$ of $\frak g$ containing $\frak b_{\frak k}$.  Define a section
$$s: W(\frak k, \frak t)\backslash W^{\theta}(\frak g, \frak h) \rightarrow W^{\theta}(\frak g, \frak h)$$ by $s(W(\frak k,\frak t)\cdot w)=w_{\frak k}w$, where $w_{\frak k}$ is the unique element on $W(\frak k, \frak t)$ such that $w_{\frak k} w(\frak b) \cap \frak k = \frak b_{\frak k}$. Then the image of $s$ is $ W^{\theta}_{1}  :=\{ w \in W^{\theta}(\frak g, \frak h)| w(\frak b) \cap \frak k  = \frak b_{\frak k} \}$.

%Let $\frak b_{\frak k}$ be a Borel subalgebra of $\frak k$.

\begin{proposition} [Section 18.3 of \cite{FHW}]
The space of base cycles in $G/B$ is parameterized by $W(\frak k, \frak t)\backslash W^{\theta}(\frak g, \frak h) \simeq W^{\theta}_{1}$.

\end{proposition}

 Once we fix a Borel subalgebra $\frak b_{\frak k}^{\rm ref}$ of $\frak k$ and a Borel subalgebra $\frak b^{\rm ref}$ of $\frak g$ containing $\frak b_{\frak k}^{\rm ref}$,
 the $K$-orbits $C(w)=K\cdot[wB^{\rm ref}]$, $w \in W_1^{\theta} $, are all base cycles in $G/B$, and $\pi(C(w)) $, $w \in W_1^{\theta}  $, are all base cycles in $G/Q$, where $\pi:G/B \rightarrow G/Q$ is the projection map.

\subsection{Extremal weights and their indices} \label{sect:extremal weights and indices}

\subsubsection{Hasse diagram} \label{Hasse diagram}

 Let $W:=W(\frak k , \frak t)$ and let $\omega_0$ be the longest element of $W$. Let $\{\beta_1, \dots, \beta_{s}\}$ be the simple root system of $(\frak k, \frak b_{\frak k}^{\rm ref})$ and $\{\xi_1, \dots, \xi_s\}$ be the fundamental weight system. For a subset $I$ of $\{ 1, \dots, s\}$, let $\lambda_I = \sum_{i \not\in I} \xi_i$ and let $\frak p_I=\frak l_I + \frak u_I$ be the Levi decomposition of the parabolic subalgebra of $\frak k$ generated by positive roots of $\frak k$ and $\beta_i, i \in I$.  For $w \in W$, put $\Delta(w)=\{ \alpha \in \Delta^+: w^{-1}\alpha  \in - \Delta^+\}$.  Then the orbit $W.\lambda_I$ of $\lambda_I$ is identified with $W^I:=\{w \in W: \Delta(w) \subset \Delta(\frak u_I)\}$, a (directed) subgraph of $W$,  called the Hasse diagram attached to the standard parabolic subalgebra $\frak p_I$ of $\frak k$.

 We remark that $W^I$ parameterizes $B_{\frak k}$-orbits in $K$-flag manifold $K/P_I$, {and also parametrizes} the highest weight orbit of the irreducible $K$-representation of highest weight $\lambda_I$.  For details, see Section 4.3 of \cite{BE}. 

\subsubsection{Highest weights of $\frak s$ and $\frak s_{\pm}$}    \label{extremal weights of s}

{ Let us explain the highest weight of the K-representation space $\frak s$. For details, see Section 17.4 of \cite{FHW}.  

First assume $\text{rank}\, \frak g = \text{rank}\, \frak k$. Let $\Psi = \{\psi_1,\ldots, \psi_r\}$ be the simple root system of $ (\frak g, \frak h )$  and let $\psi_0$ be the negative of the maximal root.
Write $-\psi_0=\sum n_i\psi_i$.
If $G_0$ is of Hermitian type, then $\frak k = \frak z\oplus \frak k'$, where $\frak k'=[\frak k, \frak k]$ is the semisimple component of $\frak k$, $\frak z$ is the center with $\dim \frak z=1$. The tangent space representation is {the} direct sum $\nu\oplus\nu^*$ of {the} irreducible representation $\nu$ and its dual $\nu^*$. 
%, where $\nu$ is the representation on $\frak s_+$. 
%The representation $\nu$ is irreducible, 
The highest weight of the representation $\nu$ is the maximal root $-\psi_0$, and $\nu = \xi + \nu'$ where $\xi = -\psi_0|_{\frak z}$ and $\nu'$ represents $\frak k'$.
If $G_0$ is of non-Hermitian type, $\frak k$ is semisimple. The tangent space representation is the irreducible representation {$\nu_{-\psi_j}$} of $\frak k$ of highest weight $-\psi_j$ with $n_j=2$ (Theorem 17.4.4 of \cite{FHW}).

Assume $\text{rank}\, \frak g > \text{rank}\, \frak k$. Then $\frak k$ is semisimple and the tangent space representation of $\frak k$ on $\frak s$ is the adjoint representation which is irreducible. (Theorem 17.4.7 of \cite{FHW}).} \\

In addition, we have the following properties. \\

$\bullet$ When $G_0$ is not of Hermitian type, $\frak s$ is an irreducible $K$-representation space. The highest weight $\lambda_{\frak s}$ of $\frak s$ is given {in terms of fundamental weights} in Table 17.4.6 and Table 17.4.8 of \cite{FHW}. 
Moreover, $-w_0\lambda_{\frak s}$ is the highest weight of $\frak s^*$ with the longest element $w_0\in W$ and $\frak s^* = \frak s$.
From this list we see that  (1) $\frak s$ is self-dual, i.e., $ \omega_0\lambda_{\frak s} = -\lambda_{\frak s}$, and (2) $\lambda_{\frak s}$ are of the form $\lambda_I$ or its multiple or their sums.  

{ Note that the elements of the $i$-th column of  the matrix $(A^t)^{-1}$  inverse  to the transposed Cartan matrix $A$   are the coefficients of simple roots when we express the $i$-th fundamental weight as a linear combination of simple roots. We use the matrix $(A^t)^{-1}$   to  get a  description of $\lambda_{\frak s}$  as a  linear combination of simple roots of $\frak k$.}   
  \\

$\bullet$ When $G_0$ is of Hermitian type,
 the action of $\frak k$ on $\frak s$ decomposes into the direct sum $\frak s =\frak s_-+ \frak s_+$ of two irreducible subrepresentations.
The highest weight $\lambda_+$ of $\frak s_+$ is given {in terms of fundamental weights} in Table 17.4.5 of \cite{FHW}. Since $\frak s_-$ is the dual representation space of $\frak s_+$, the highest weight $\lambda_-$ of $\frak s_-$ is $-\omega_0 \lambda_+$. Thus $\omega_0\lambda_+ = -\lambda_-$.  As in non-Hermitian cases, $\lambda_{\pm}$ are of the form $\lambda_I$ or its multiple or their sums. \\

%\subsubsection{Various Weyl groups} \label{various Weyl groups}

  \subsubsection{Fibers $E_0^{w_1}$ at base points of base cycles} \label{index of E0}

   Fix a Borel subalgebra $\frak b_{\frak k}^{\rm ref}$ of $\frak k$ and a Borel subalgebra $\frak b^{\rm ref}$ of $\frak g$ containing $\frak b_{\frak k}^{\rm ref}$. 
 
  % We will use simple root systems $\Psi_{\frak k}$ of $\frak k$ and $\Psi$ of $\frak g$ and its restriction $\Psi_{\theta}$ to $\frak t$, the subset $W_1^{\theta}$ of the Weyl group $W(\frak g, \frak h)$ described in  Chapters 19--20 of \cite{FHW}. For $\frak g_0$ which are not dealt with in Chapters 19--20 of \cite{FHW} we will describe them below. 

At the base point of $[B^{\rm ref}] \in X=G/B$, the central fiber $E_0$ of the normal bundle $E$ of the base cycle $K\cdot[B^{\rm ref}]$ is
$$E_0 = \frak s/(\frak s \cap \frak b^{\rm ref})=\frak s \cap \frak b^{{\rm ref}, n}.$$  
Here, $\frak b^{\text{ref},n}$ is the nilradical of $\frak b^{\text{ref}}$.
For $w_1 \in W_1^{\theta}$,
 the fiber $E_0^{w_1}$ at $[w_1(B^{\rm ref})] \in X=G/B$ is
 $$E_0^{w_1}=\frak s /\frak s \cap w_1(\frak b^{\rm ref})=\frak s \cap w_1(\frak b^{{\rm ref},n}).$$

 As in Section \ref{sect:classification of base cycles}, choose   a $\tau$-stable $\theta$-stable  Cartan subalgebra $\frak h$ of $\frak g$ such that $\frak h_0:=\frak h \cap \frak g_0$ is a maximally compact Cartan subalgebra of $\frak g_0$, i.e., $\frak t:=\frak h_{\frak k}=\frak h \cap \frak k$ is a Cartan subalgebra of $\frak k$. Let $\Sigma(\frak g, \frak h)$ denote the set of roots of $\frak g$ with respect to $\frak h$. By restriction, we get a map from $\Sigma(\frak g, \frak h)$ to  the set $M(\frak t, \frak g) \backslash \{0\}$ of nonzero $\frak t$-weights of the representation $\frak g$ of $\frak k$, which is again a root system $\Sigma(\frak g, \frak t)$ (Proposition 18.2.3 of \cite{FHW}).  
 Furthermore, $W^{\theta}(\frak g, \frak h)$,   the subgroup of $W(\frak g, \frak h)$ consisting of $w \in W(\frak g, \frak h)$ such that $w \circ \theta = \theta \circ w$, acts on $\Sigma(\frak g, \frak t)$. Therefore, $W_1^{\theta}$ also acts on  $\Sigma(\frak g, \frak t)$.  

{  Each element $w_1 \in W_1^{\theta}$ determines  
 positive root systems $\Sigma^{+,w_1}(\frak g, \frak h)$, $\Sigma^{+,w_1}(\frak g, \frak t)$, $\Sigma^{+,w_1}(\frak k, \frak t)$ consistent with the restriction map 
 ${\rm res}: \Sigma(\frak g, \frak h) \rightarrow \Sigma(\frak g, \frak t)$ and the inclusion map  $\Sigma(\frak k,\frak t) \hookrightarrow \Sigma(\frak g, \frak t) $, so that we have 
 $$ {\rm res}: \Sigma^{+, w_1}(\frak g, \frak h)  \rightarrow \Sigma^{+, w_1}(\frak g, \frak t) \text{ and } \Sigma^{+,w_1}(\frak k, \frak t) \hookrightarrow\Sigma^{+,w_1}(\frak g, \frak t).$$}
 %the intersection of $\Sigma(\frak g, \frak t)$ with the set $M(\frak t, w_1(\frak b^{{\rm ref},n}))$ of $\frak t$-weights in $w_1(\frak b^{{\rm ref} ,n})$. 
 In particular, the set of $\frak t$-weights in $E_0^{w_1}$  is a subset of the positive roots system $\Sigma^{+,w_1}(\frak g, \frak t)$.

\medskip 

Recall that, for two weights $\lambda$ and $\lambda'$, by the distance from $ \lambda $ to $ \lambda'$, we mean the minimum of $\ell(w)$ with $w \lambda  =  \lambda'$ and $w \in W$ (Section \ref{sect: ampleness formula}).  

%$\bullet$ $E_0 = \frak s/(\frak s \cap (\frak q + \theta (\frak q))$ the fiber at $Q \supset  B^{\rm ref}$.

%$E_0 = \frak s \cap \frak q^n$ at $Q  \in G/Q$ such that  $Q \supset B^{\rm ref}$ if measurable.

\begin{proposition}  \label{prop: strategy} $\,$ \\
\begin{enumerate}
\item When $G_0$ is not of Hermitian type, for $w_1 \in W_1^{\theta}$,
 % let $w_{\max}$ be an element in $W$ of maximal length among all $w \in W$ such that $w_1^{-1}(w \lambda_{\frak s}) $ is a positive root with respect to $ \frak b^{{\rm ref},n}$. 
  %
  %Let 
  let $w_{\max}$ be an element in $W$  such that the maximum distance between $\lambda_{\frak s}$ and $w \lambda_{\frak s} \in \Lambda(w_1(\frak b^{{\rm ref},n}))$ where $w \in W$,  is attained at $w_{\max}\lambda_{\frak s}$.  
  Then $$ \ind (-\Lambda_{\ext}(E_0^{w_1}))=\ind (-w_{\max}\lambda_{\frak s}) . $$

\item  When $G_0$ is of Hermitian type,  for $w_1 \in W_1^{\theta}$, 
let $w_{\max}$ be an element in $W$  such that the maximum distance between $\lambda_{\pm}$ and $w \lambda_{\pm} \in \Lambda(w_1(\frak b^{{\rm ref},n}))$ where $w \in W$,  is attained at $w_{\max}\lambda_{\pm}$.  
%
%let   $w_{\max, \pm}$ be an element in $W $ of maximal length  among all    $w \in W$ such that  $w_1^{-1}(w\lambda_{\pm})  $ is a positive root with respect to $\frak b^{{\rm ref},n}$. 
Then $$ \ind (-\Lambda_{\ext}(E_0^{w_1}))=\min\{ \ind (-w_{\max,+}\lambda_{+}), \ind (-w_{\max,-}\lambda_{-})\}.$$
In particular, for $w_1=id$, we have $ \ind (-\Lambda_{\ext}(E_0))=0$.
\end{enumerate}

\end{proposition}
\begin{proof}
(1) When $G_0$ is not of Hermitian type, we have 
 $\Lambda_{\max}(E_0)=\{\lambda_{\frak s}\}$,   so that  
  $$\Lambda_{\ext} (E_0) = \{w\lambda_{\frak s}: w \in W, w\lambda_{\frak s} \in \Lambda(E_0)\} .$$    Among weights in $W.\lambda_{\frak s}$,   a unique dominant one is $\lambda_{\frak s}$. Thus, for a weight $-\lambda  $  in $-W.\lambda_{\frak s}$, $w(-\lambda) $ is dominant if and only if $w(-\lambda)=\lambda_{\frak s}$, or equivalently, $w \lambda =-\lambda_{\frak s}=w_0\lambda_{\frak s}$ because $\frak s$ is self-dual. 
  Therefore,  the minimum $\ind (-\Lambda_{\ext}(E_0))$ occurs at $-\lambda \in - W.\lambda_{\frak s} $  
such that the distance between $\lambda$ and $w_0\lambda_{\frak s}$ is minimal among $\lambda \in \Lambda(E_0)$, in other words, the distance between $\lambda_{\frak s}$ and $\lambda$ is maximal among $\lambda \in \Lambda(E_0)$.     
  
  Let $w_{\max} $ be an element in $W $ 
    such that 
  the maximum   distance between $\lambda_{\frak s}$ and $u\lambda_{\frak s} \in \Lambda(E_0)$ where $u \in W$, is attained at $w_{\max}\lambda_{\frak s}$. 
  %$w_{\max}\lambda_{\frak s}$ is still in $\Lambda(E_0)$ and the distance between $\lambda_{\frak s}$ and $w_{\max}\lambda_{\frak s}$ is maximal among all such $\lambda$.} 
  %of maximal length  among all $u \in W$ such that  $u\lambda_{\frak s}  $ is still in {$\Lambda(E_0)$}. %, i.e., is a positive root with respect to $\frak b^{{\rm ref}, n}$. %{ Recall that the index $\ind(\lambda)$ of a weight $\lambda$ is the minimum length $\ell(w)$ of $w \in W$ such that $w\lambda$ is dominant.  
 % { Thus} the minimum $\ind (-\Lambda_{\ext}(E_0))$ occurs at $-w_{\max}\lambda_{\frak s}$. Thus
 Then   by the arguments in the above, we get  $ \ind (-\Lambda_{\ext}(E_0))=\ind (-w_{\max}\lambda_{\frak s})$.   % Let $K/P_I$ be the highest  weight orbit of  $\frak s$.
% By   symmetry (or/and self-duality $w_0\lambda_{\frak s} =-\lambda_{\frak s} $),
% \begin{eqnarray*}
 % \ind (-\Lambda_{\ext}(E_0))= \left\{
 % \begin{array}{cc}
 % r &\text{ if } \dim K/P_I=2r-1 \\
 % r-1 &\text{ if } \dim K/P_I=2r-2.
 % \end{array}\right.
%\end{eqnarray*}
%if $\frak k$ is simple, is the minimum if $\frak k$ is semisimple (why?) \\

If $w_1 \not=id \in W_1^{\theta}$, we have $\Lambda_{\max}(E_0^{w_1})=\{\lambda_{\frak s}\}$ and hence  
 %  (why? maybe because $\frak b_{\frak k}^{\rm ref} \subset w_1(\frak b^{\rm ref})$). 
$$\Lambda_{\ext}(E_0^{w_1}) =W.\lambda_{\frak s} \cap \Lambda(\frak s \cap w_1(\frak b^{{\rm ref},n}))=W.\lambda_{\frak s} \cap  \Lambda(w_1(\frak b^{{\rm ref},n}))$$
 because $W.\lambda_{\frak s} \subset \Lambda(\frak s)$. 
%  { Furthermore, 
%(why?
%   $$\max\{\ell(w): w\lambda_{\frak s} \in \Lambda(E_0^{w_{1}}), w \in W\} = \max\{\ell(w):w_{1}^{-1}(w\lambda_{\frak s}) \in \Lambda(\frak b^{{\rm ref},n}), w \in W\}.$$}
  %) 
   Let $w_{\max}$ be an element in $W$  such that the maximum distance between $\lambda_{\frak s}$ and $w \lambda_{\frak s} \in \Lambda(w_1(\frak b^{{\rm ref},n}))$ where $w \in W$,  is attained at $w_{\max}\lambda_{\frak s}$.  
  Then $ \ind (-\Lambda_{\ext}(E_0^{w_1}))=\ind (-w_{\max}\lambda_{\frak s})$.  
   
%\\

(2) When $G_0$ is of Hermitian type, $\Lambda_{\max}(E_0) =\{\lambda_+\}$.  
Let $w_{\max}$ be an element in $W$  such that the maximum distance between {$\lambda_{+}$} and $w \lambda_{+} \in \Lambda( \frak b^{{\rm ref},n}) $ where $w \in W$,  is attained at $w_{\max}\lambda_{+}$.  
%Let $w_{\max}$ be an element in $W $ such that  the maximum distance between $\lambda_+$ and $w\lambda_{+}  $ where $w \in W$ and $w \lambda_+$  is positive with respect to $\frak b^{{\rm ref},n}$, is attained at $w_{\max}\lambda_+$. 
Then $ \ind (-\Lambda_{\ext}(E_0))=\ind (-w_{\max}\lambda_{+})$. Since $w_{\max}\lambda_+=-\lambda_-$, we obtain $ \ind (-\Lambda_{\ext}(E_0))=0$.   Note that in this case $w_{\max}=w_0$.

If $w_1 \not=id \in W_1^{\theta}$, then $\Lambda_{\max}(E_0^{w_1})\subset \{\lambda_{+}, \lambda_-\}$. %(why?). 
Let {$w_{\max,\pm}$} be an element in $W$  such that the maximum distance between $\lambda_{\pm}$ and $w \lambda_{\pm} \in \Lambda(w_1(\frak b^{{\rm ref},n}))$ where $w \in W$,  is attained at {$w_{\max,\pm}\lambda_{\pm}$}.  
%Let $w_{\max, \pm}$ be an element in $W $ of maximal length  among all $w \in W$ such that  $w_1^{-1}(w\lambda_{\pm})  $ is positive with respect to $\frak b^{{\rm ref},n}$. 
Then $ \ind (-\Lambda_{\ext}(E_0^{w_1}))=\min\{ \ind (-w_{\max,+}\lambda_{+}), \ind (-w_{\max,-}\lambda_{-})\}$.   %(why?) \\
\end{proof}

We remark that the condition $w \lambda_{\frak s} \in \Lambda(w_1(\frak b^{{\rm ref},n}))$ is equivalent to that  $w_1^{-1}(w \lambda_{\frak s}) $ is a positive root with respect to $ \frak b^{{\rm ref},n}$ and the condition $w \lambda_{\pm} \in \Lambda(w_1(\frak b^{{\rm ref},n}))$ is equivalent to that  $w_1^{-1}(w \lambda_{\pm}) $ is a positive root with respect to $ \frak b^{{\rm ref},n}$. 
 In the computation in Section  \ref{sect:full flag manifolds}, we will  use these conditions to find out $w_{\max}\lambda_{\frak s} $ and $w_{\max}\lambda_{\pm}$.  
 
\subsection{Our strategy} \label{sect:our strategy}

We will  compute $\dim C^{w_1}-a(E^{w_1})=  \ind (-\Lambda_{\ext}(E_0^{w_1}))$, where $E^{w_1}$ is the normal bundle  of the base cycle $C^{w_1}:=K\cdot[w_1(B^{\rm ref})]$ in $X=G/B$, for $w_1 \in W_1^{\theta}$. \\

 Fix a Borel subalgebra $\frak b_{\frak k}^{\rm ref}$ of $\frak k$ and a Borel subalgebra $\frak b^{\rm ref}$ of $\frak g$ containing $\frak b_{\frak k}^{\rm ref}$ in a standard way as in Chapters 19--20 of \cite{FHW}. We will use simple root systems $\Psi_{\frak k}$ of $\frak k$ and $\Psi$ of $\frak g$, % and its restriction $\Psi_{\theta}$ to $\frak t$, 
 the subset $W_1^{\theta}$ of the Weyl group $W(\frak g, \frak h)$ described in  Chapters 19--20 of \cite{FHW}. For $\frak g_0$ which are not dealt with in Chapters 19--20 of \cite{FHW} we will describe them below.  \\

\begin{enumerate}
\item [\bf Step 1.] Describe $W.\lambda_{\frak s}$ or $W.\lambda_+$ as a (directed) subgraph $W^I$ of $W$ (Section \ref{Hasse diagram}).

\item [\bf Step 2.] From   Table 17.4.6 and Table 17.4.8 of \cite{FHW} and 
Table 17.4.5 of \cite{FHW}, 
we get an expression of $\lambda_{\frak s}$  as a restricted root of $\frak g$.  %$\frak g$ 
Here, we use the matrix $(A^t)^{-1}$  inverse  to the transposed Cartan matrix $A$ ({Table 2 of} \cite{OV})  and the root system $\Psi_{\frak k}$ of $\frak k$  %$\Psi$ of $\frak g$ 
(Chapter 19--20 of \cite{FHW}).  
(Section \ref{extremal weights of s})  
Applying $w \in W^I$ obtained in {\bf Step 1.} to $\lambda_{\frak s}$, we get an expression of $w\lambda_{\frak s}$ or $w\lambda_+$, where $w \in W$, as a restricted root of $\frak g$. 

\item[\bf Step 3.] For $w_1 \in W_1^{\theta}$, describe $w_1^{-1}(w\lambda_{\frak s})$  or $w_1^{-1}(w\lambda_{\pm})$, where $w \in W$, as a {restricted} root in $\frak g$.
    (The action of $W_1^{\theta}$ on {restricted} roots of $\frak g$ is described in Chapter 19--20 of \cite{FHW}.)  

\item[\bf Step 4.]  Find $w_{\max}\lambda_{\frak s}$ or $w_{\max, \pm}\lambda_{\pm}$ (Section \ref{index of E0}) and compute  $\ind(-w_{\max}\lambda_{\frak s})$ or $\ind (-w_{\max,+}\lambda_{+})$, $ {\ind (-w_{\max,-}\lambda_{-})}$ using the (directed) subgraph $W^I$ of $W$. \\

\end{enumerate}

\section{Ampleness of base cycles on full flag varieties} \label{sect:full flag manifolds}

From what follows $M_{n}(\mathbb C)$ denotes
the set of $n \times n$ matrices with complex coefficients and $B_n^-$ denotes the set of $n\times n$ lower triangular matrices.

 \subsection{$\frak g_0=\frak{sl}(m,\mathbb R)$}\label{sl}
{The Lie algebra $\frak g = \frak{sl}(m,\mathbb R)$ and $\frak b^{\text{ref}} = \frak{sl}(m,\mathbb C)\cap B_m^-$.
The associated flag manifold $G/B$ consists of all full flags
$$
0\subset F_1\subset F_2\subset\cdots\subset F_{m-1} \subset F_m = \mathbb C^m
\quad\text{ with }
\dim_{\mathbb C} F_j = j \text{ for all } j=1,\ldots, m.
$$
For a complex conjugate $\tau\colon \mathbb C^m\to\mathbb C^m$ defined by $\tau(z) = \overline z$, a flag $(0\subset F_1\subset\cdots\subset F_{m-1}\subset\mathbb C^m)$ is said to be $\tau$-generic, if $\dim (F_i\cap \tau(F_j)) = \max\{0, i+j-m\}$.
In \cite{HSB01}, it was proved that if $m$ is an odd number, then the unique open $SL(m,\mathbb R)$-orbit is the set of $\tau$-generic flags, and if $m$ is an even number, then depending on the orientation on $\mathbb C^m_{\mathbb R}$, there are two open orbits which are the set of positively oriented $\tau$-generic flags and the set of negatively oriented $\tau$-generic flags.}

Let $b(z,w) := z^t R_m w$ define a nondegenerate symmetric bilinear form on $\mathbb C^m$ where $R_m$ denotes the $m\times m$ symmetric matrix with $1$ on the antidiagonal and $0$ elsewhere.
The Lie algebra $\frak g=\frak{sl}(m, \mathbb C)$ is the space of trace zero matrices in $M_{m}(\mathbb C)$ and $\frak k = \frak{so}(m,\mathbb C)$ is the Lie algebra of the isotropy group of $b$, which is the set of fixed points under the involution $\theta (X) = R_m(-X^t) R_m$.
Let $\frak h :=\left\{ \text{diag}(\delta_1,\ldots, \delta_m)\in M_m(\mathbb C) : \delta_1+\cdots + \delta_m=0\right\}$ be a Cartan subalgebra of $\frak{sl}(m,\mathbb R)$
and $\frak t:= \frak h\cap \frak{so}(m,\mathbb C)$ be the Cartan subalgebra of $\frak k$. The reference Borel subalgebras are
$\frak b_{\frak k}^{\text{ref}} = \frak{so}(m,\mathbb C)\cap B_m^-$.
{%The number of $SL(m,\mathbb R)$-flag domains in $G/B$ is one if $m$ is even and two if $m$ is odd (cf. \cite{FHW} Chapter 13). 
In all cases of $\frak g_0 = \frak{sl}(m,\mathbb R)$, we have $\text{rank}(\frak g)>\text{rank}(\frak k)$ and $\lambda_{\frak s}$ is the restriction of the longest root of $\Sigma(\frak g, \frak h)$.}

\medskip

 \noindent \textsf{\textbf{Case 1.}} If $m=2r+1$ with $r \geq 2$,
$\frak t$ is given by $\{\text{diag}(\ve_1,\ldots, \ve_r,0,-\ve_r, \ldots, -\ve_1): \ve_j\in \mathbb C \text{ for all } j\}$.
With respect to $\frak t$, the simple root system of $\frak k$ is given by
$$
\{ \ve_1 - \ve_2,
\ve_2 - \ve_3, \ldots,
\ve_{r-1} - \ve_r, \ve_r\}
=: \{ \beta_1,\beta_2,\ldots, \beta_{r-1}, \beta_r\}
$$
and the fundamental weight system is given by
$\{\xi_1, \xi_2,\ldots, \xi_r\}$ where
\begin{equation}\nonumber
\begin{aligned}
\xi_i &= \beta_1+2\beta_2+\cdots+ (i-1)\beta_{i-1} + i(\beta_i+\beta_{i+1}+\cdots + \beta_r) \text{ with } i<r, \text{ and } \\
\xi_r &= \frac{1}{2}(\beta_1 + 2\beta_2 +\cdots+r\beta_r).
\end{aligned}
\end{equation}
Note that for $1\leq i \leq r$,
$\ve_i=\beta_i+\beta_{i+1} + \cdots + \beta_r$.
Since $\lambda_{\frak s} =2 \xi_1=2 \ve_1$, the Hasse diagram of $W.\lambda_{\frak s}$ is

\begin{tikzpicture}
    \matrix (a) [matrix of math nodes, column sep=0.35cm ]{
     \lambda_{\frak s}= 2\ve_1& 2\ve_2& \cdots& 2\ve_{r-1}& 2\ve_r& -2\ve_r& -2\ve_{r-1}& \cdots& -2\ve_2& -2\ve_1=-\lambda_{\frak s}.\\
    };
    \draw[->] (a-1-1)--(a-1-2);
    \draw[->] (a-1-2)--(a-1-3);
    \draw[->] (a-1-3)--(a-1-4);
    \draw[->] (a-1-4)--(a-1-5);
    \draw[->] (a-1-5)--(a-1-6);
    \draw[->] (a-1-6)--(a-1-7);
    \draw[->] (a-1-7)--(a-1-8);
    \draw[->] (a-1-8)--(a-1-9);
    \draw[->] (a-1-9)--(a-1-10);
\end{tikzpicture}
Thus $w_{\max}\lambda_{\frak s} = 2 \ve_r$ and hence $\ind(-w_{\max}\lambda_{\frak s}) =r$.

Since we have $W_1^{\theta}=\{id\}$, i.e. there is a unique base cycle in $G/B$, $\ind(-\Lambda_{\ext}(E_0))=r$.

\medskip

 \noindent \textsf{\textbf{Case 2.}}  If $m=2r$ with $r \geq 3$,
$\frak t$ is given by $\{\text{diag}(\ve_1,\ldots, \ve_r,-\ve_r, \ldots, -\ve_1): \ve_j\in \mathbb C \text{ for all } j\}$.
The simple root system of $(\frak t, \frak k)$ is given by
$$
\left\{ \ve_1 - \ve_2,
\ve_2 - \ve_3, \ldots,
\begin{array}{c}
\ve_{r-1} - \ve_r\\
\ve_{r-1}+\ve_r
\end{array}
\right\}
=: \left\{ \beta_1,\beta_2,\ldots,
 \begin{array}{c}
\beta_{r-1}\\
\beta_r
\end{array}
\right\},
$$
the simple system for $(\frak t, \frak g)$ is given by
$$
\{\ve_1-\ve_2, \ve_2-\ve_3, \ldots, \ve_{r-1}-\ve_r, 2\ve_r\} =:\{ \gamma_1, \gamma_2, \ldots, \gamma_r\}
$$
and the fundamental weight system is given by
$\{\xi_1, \xi_2,\ldots, \xi_r\}$ where
\begin{equation}\nonumber
\begin{aligned}
\xi_i &= \beta_1+2\beta_2+\cdots+ (i-1)\beta_{i-1} + i(\beta_i+\cdots + \beta_{r-2}) + \frac{1}{2} i(\beta_{r-1}+\beta_r) \text{ with } i\leq r-2,  \\
\xi_{r-1} &= \frac{1}{2}\left(\beta_1 + 2\beta_2 +\cdots+(r-2)\beta_{r-2}+\frac{1}{2}r\beta_{r-1} + \frac{1}{2}(r-2)\beta_r
\right), \text{ and }\\
\xi_{r} &= \frac{1}{2}\left(\beta_1 + 2\beta_2 +\cdots+(r-2)\beta_{r-2}+\frac{1}{2}(r-2)\beta_{r-1} + \frac{1}{2}r\beta_r
\right).
\end{aligned}
\end{equation}
Note that for $1 \leq i \leq r-2$, $2 \ve_i =2 \beta_i + \cdots + 2 \beta_{r-2} + \beta_{r-1} + \beta_r$,
$2 \ve_{r-1} =\beta_{r-1} + \beta_r$, and
$2 \ve_r=-\beta_{r-1} + \beta_r=\gamma_r$.
Since $\lambda_{\frak s} =2 \xi_1=2 \ve_1$ and the Hasse diagram of $W.\lambda_{\frak s}$ is

\begin{tikzpicture}
    \matrix (a) [matrix of math nodes, column sep=0.35cm ]{
    &&&&2\ve_r&&&\\
     \lambda_{\frak s}= 2\ve_1
     & 2\ve_2
     & \cdots
     & 2\ve_{r-1}
     &~
     & -2\ve_{r-1}
     & \cdots
     & -2\ve_2
     & -2\ve_1=-\lambda_{\frak s}.\\
     &&&&-2\ve_r&&&\\
    };
    \draw[->] (a-2-1)--(a-2-2);
    \draw[->] (a-2-2)--(a-2-3);
    \draw[->] (a-2-3)--(a-2-4);
    \draw[->] (a-2-4)--(a-1-5);
    \draw[->] (a-2-4)--(a-3-5);
    \draw[->] (a-1-5)--(a-2-6);
    \draw[->] (a-3-5)--(a-2-6);
    \draw[->] (a-2-6)--(a-2-7);
    \draw[->] (a-2-7)--(a-2-8);
    \draw[->] (a-2-8)--(a-2-9);
\end{tikzpicture}
 
\noindent Thus $w_{\max}\lambda_{\frak s} = 2 \ve_r$ and $\ind(-w_{\max}\lambda_{\frak s}) =r-1$. \\

% Note that the highest weight orbit in $\mathbb P(\frak s)$ is $Q^{2r-2}$ in the second Veronese embedding.

 Since $W_1^{\theta}=\{id, s_{\gamma_r}\}$, i.e., there are two base cycles in $G/B$, we need to consider the normal bundle of the base cycle $C^{s_{\gamma_r}} = K\cdot [s_{\gamma_r}(B^{\text{ref}})]$.
The simple root systems associated to $s_{\gamma_r}(\frak b^{\rm ref})$ is given by
\begin{equation}\nonumber
\begin{aligned}
& \tilde \gamma_1=\gamma_1, \ldots ,
\tilde \gamma_{r-2} = \gamma_{r-2},\, 
\tilde \gamma_{r-1} =  \gamma_{r-1} + \gamma_r, \,
\tilde \gamma_r = -\gamma_r\\
& \tilde \beta_1=\beta_1, \ldots , 
\tilde \beta_{r-2} = \beta_{r-2},\,
\tilde \beta_{r-1} =  \beta_{r},\,
\tilde \beta_r = \beta_{r-1}
\end{aligned}
\end{equation}
and   $W.\lambda_{\frak s}$ is

\begin{tikzpicture}
    \matrix (a) [matrix of math nodes, column sep=0.35cm ]{
    &&&&2\ve_r&&&\\
     \lambda_{\frak s}= 2\ve_1
     & 2\ve_2
     & \cdots
     & 2\ve_{r-1}
     &~
     & -2\ve_{r-1}
     & \cdots
     & -2\ve_2
     & -2\ve_1=-\lambda_{\frak s}.\\
     &&&&-2\ve_r&&&\\
    };
    \draw[->] (a-2-1)--(a-2-2);
    \draw[->] (a-2-2)--(a-2-3);
    \draw[->] (a-2-3)--(a-2-4);
    \draw[->] (a-2-4)--(a-1-5);
    \draw[->] (a-2-4)--(a-3-5);
    \draw[->] (a-1-5)--(a-2-6);
    \draw[->] (a-3-5)--(a-2-6);
    \draw[->] (a-2-6)--(a-2-7);
    \draw[->] (a-2-7)--(a-2-8);
    \draw[->] (a-2-8)--(a-2-9);
\end{tikzpicture}
 
\noindent Since 
\begin{equation}\nonumber
\begin{aligned} 
&2 \ve_i =2 \tilde \beta_i + \cdots +   2\tilde\beta_{r-2} +  \tilde\beta_{r-1} +  \tilde\beta_r \textup{ for } 1 \leq i \leq r-2\\
&2 \ve_{r-1} =\tilde\beta_{r-1} + \tilde\beta_r, \quad \text{  and }\\
&2 \ve_r = -2\tilde \beta_{r-1}+  2 \tilde\beta_r= -\gamma_r,
\end{aligned}
\end{equation}
we obtain $w_{\max}\lambda_{\frak s} =- 2 \ve_r$ and hence $\ind(-w_{\max}\lambda_{\frak s}) =r-1$.
As a result $\ind(-\Lambda_{\ext}(E_0^{w_1}))=r-1$ for any $w_1 \in W_1^{\theta}$.

\medskip

 \noindent \textsf{\textbf{Case 3.}} If $m=3$, {then $\frak k=\frak{so}(3,\mathbb C)$ and $\frak t$ is given by 
 $\{\text{diag}(\ve_1,0,-\ve_1) : \ve_1\in \mathbb C\}$ as Case~1. The simple root system of $(\frak t, \frak k)$ is $\Psi_{\frak k} = \{\ve_1\}$ and {$\Sigma(\frak g, \frak t) = \{\pm\ve_1, \pm 2\ve_1\}$}.
Since $\lambda_{\frak s} = 2\xi_1=2\ve_1$, the Hasse diagram of $W.\lambda_{\frak s}$ is 

\begin{center}
\begin{tikzpicture}
    \matrix (a) [matrix of math nodes, column sep=0.35cm ]{
    \lambda_{\frak s} = 2\ve_1& -2\ve_1 = -\lambda_{\frak s}\\
    };
    \draw[->] (a-1-1)--(a-1-2);
\end{tikzpicture}
\end{center}
}

\noindent
Since $w_\text{max} \lambda_{\frak s} = 2\ve_1$, we have
$\text{ind}(-w_\text{max}\lambda_{\frak s} )= 1$
and hence $\text{ind}(-\Lambda_\text{ext}(E^{w_1}_0)) = 1$ for any $w_1\in W^\theta_1$.
Note that $W_1^{\theta}=\{id\}$.
{We remark that this case is essentially the same as Case~1 except for the type of the Dynkin diagram.}

\medskip

\noindent \textsf{\textbf{Case 4.}} If $m=4$, {$\frak{sl}(4,\mathbb R)\cong \frak{so}(3,3)$. Thus, we will consider this case in the section~\ref{so(m,n)} Case 1. }

\medskip

\noindent \textsf{\textbf{Case 5.}} If $m=2$, the case is trivial since the only relevant flag manifold is the projective line $\mathbb C\mathbb P^1$,
the cycles are points, and flag domains are biholomorphic
to the unit disc in complex plane.

\medskip

As a result, for $\frak g_0 = \frak{sl}(m,\mathbb R)$ we have $\text{ind}(-\Lambda_{\text{ext}}(E_0^{w_1})) = \left[\frac{m-1}{2}\right]$.

 \subsection{$\frak g_0 =\frak{su}(p,q)$}
{In this case $G_0$ is of Hermitian type.}
{We have $\frak g = \frak{sl}(r,\mathbb C)$ with $r:=p+q$ and $\frak b^{\text{ref}}=\frak{sl}(r,\mathbb C)\cap B_r^-$. Let $I_{p,q}$ be the standard Hermitian form of signature $(p,q)$ defined by the matrix $$
\left(\begin{array}{cc}
-I_q&0\\
0&I_p
\end{array}\right).
$$
The associated flag manifold $G/B$ consists of all full flags
$$
0\subset F_1\subset F_2\subset\cdots\subset F_{r-1} \subset F_r = \mathbb C^r
\quad\text{ with }
\dim_{\mathbb C} F_j = j \text{ for all } j=1,\ldots, r.
$$
Let $a$ be a sequence $0\leq a_1\leq\cdots\leq a_r\leq q$ and $b$ be a sequence ${0}\leq b_1\leq \cdots \leq b_r\leq p$ such that $a_j+b_j=j$ for each $j=1,\ldots, r$. 
For a fixed pair $(a,b)$ the open orbit $D_{a,b}$ in $G/B$ is described by the set of all full flags $(0\subset F_1\subset\cdots\subset F_{r-1}\subset \mathbb C^r)$ such that $F_j$ has signature $(a_j,b_j)$ with respect to $I_{p,q}$.
There are $\frac{(p+q)!}{p!q!}$ number of flag domains in $G/B$ and each orbit is determined by the signature (cf. Chapter 5.5A in \cite{FHW}).
}

Since we have $\frak k = \frak{s}(\frak{gl}(p,\mathbb C)\times\frak{gl}(q,\mathbb C))$, the Cartan subalgebras are given by
$\frak h = \frak t = \{\text{diag}(\ve_1,\ldots, \ve_r) : \ve_1+\cdots+\ve_r=0\}$.
The reference Borel subalgebra $\frak b_{\frak k}^{\text{ref}}$ is given by  $\frak{s}(\frak{gl}(p,\mathbb C)\times\frak{gl}(q,\mathbb C))\cap B^-_r$. $\lambda_+$ is the longest root of $\Sigma(\frak g, \frak h)$.
The simple root system of $\frak k$ is given by $\Psi_{\frak k}^{(1)} \cup \Psi_{\frak k}^{(2)}$ where
\begin{equation}\nonumber
\begin{aligned}
\Psi_{\frak k}^{(1)} &= \{ \ve_1-\ve_2, \ldots, \ve_{p-1}-\ve_p\}=: \{ \beta_1^{(1)},\ldots, \beta_{p-1}^{(1)}\}\quad\text{ and } \\
\Psi_{\frak k}^{(2)} &= \{ \ve_{p+1}-\ve_{p+2}, \ldots, \ve_{p+q-1}-\ve_{p+q}\}=: \{ \beta_1^{(2)},\ldots, \beta_{q-1}^{(2)}\}.
\end{aligned}
\end{equation}

 If {$1 \leq p \leq q$}, then $\lambda_+=\ve_1 - \ve_r$. {We remark that if $p=1$, then $\Psi_{\frak k}^{(1)}$ is an empty set.}
 Positive   roots of $\frak g$ are
{$\{\varepsilon_i-\varepsilon_j\}_{1\leq i<j\leq r}$}
 and the Hasse diagram of $W.\lambda_+$ is shown as Figure~\ref{su(p,q)}
 which consists of positive noncompact roots of $\frak g$.
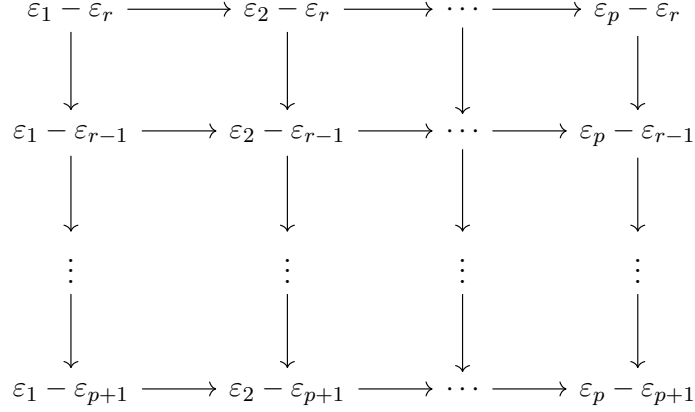
\begin{figure}[h]
 \begin{tikzcd}[row sep=1cm, column sep=1cm, nodes={anchor=center}]
\ve_1-\ve_r \rar \dar& \ve_2-\ve_r\rar \dar& \cdots \dar\rar&\ve_p-\ve_r\dar\\
\ve_1-\ve_{r-1}\rar \dar& \ve_2-\ve_{r-1}\rar \dar& \cdots\rar\dar &\ve_p-\ve_{r-1} \dar\\
    \vdots\dar&\vdots\dar&\vdots\dar&\vdots\dar\\
    \ve_1-\ve_{p+1}\rar& \ve_2-\ve_{p+1}\rar& \cdots\rar &\ve_p-\ve_{p+1}
  \end{tikzcd} 
  \caption{{Hasse diagram for} $\frak{su}(p,q)$}
  \label{su(p,q)}
\end{figure}

% \begin{eqnarray*}\begin{array}{ccccccc}
%   &   &   &\,\,\,\,\,\,\ve_{p} -\ve_{p+1}=-\lambda_- &   &        &  \\
%   &  & \ve_{p-1}-\ve_{p+1} & \ve_p - \ve_{p+2} \,\,\,  &          &  & \\
%         & \dots         & \dots & \dots \,\,\, &        &        &  \\
% \,\,\,\ve_{1}-\ve_{p+1}         & \dots         & \dots & \dots    \,\,\, &    &        &  \\
% \,\,\,\dots         & \dots         & \dots & \ve_p -\ve_r   \,\,\, &        &        &  \\
% \,\,\,\dots         & \dots         & \dots &   &        &        &  \\
% \uparrow \,\ve_1 - \ve_{r-1} &  \ve_2 - \ve_r&       &           &             &   \\
 %\,\,\,\,\lambda_+=\ve_1 - \ve_r \,\,\,\,\,\, &               &           &       &             &
% \end{array}\end{eqnarray*}

 For $a \leq p < p+1 \leq b$
 \begin{eqnarray*}
 \ind(\ve_a -\ve_b)&=& r-1-(b-a)\\
 \ind(-\ve_a + \ve_b)&=& b-a-1.
 \end{eqnarray*}
We have $w_{\max} \lambda_+ = \ve_p - \ve_{p+1}=w_0 \lambda_+$ and $-w_{\max}\lambda_+=\lambda_-$ is dominant. Thus $\ind(-w_{\max}\lambda_+) =0$.   \\

In this case, we have  $W_1^{\theta} =\{ w_{\bf j}: {\bf j} \subset \{1, \dots, r\}, |{\bf j}|=p\}$. Here, for ${\bf j}=\{j_1, \dots, j_p\} \subset \{1, \dots, r\}$ with $j_1 <j_2 < \dots <j_p$, after ordering the complement $\{1,2, \dots, r\} \backslash {\bf j}$ in increasing order, as $j_{p+1}<\dots<j_{r}$,  define $w_{\bf j} \in W^{\theta}$ by
 $$w_{\bf j}^{-1}(\ve_{ i}) =\ve_{j_i}$$
 for $1 \leq i \leq r$.   
 %Thus there are $\frac{r!}{p!q!}$ base cycles in $G/B$.
 \\

  We remark that if $j_{p+1} >j_p$, then $j_1 <j_r$. Thus we have either 
  (1)$j_1<j_r$ and $j_{p+1} >j_p$,
    or (2) $j_1>j_r$   and $j_{p+1} <j_p$, or
    (3) $j_1<j_r$   and $j_{p+1} <j_p$.  
     {{Remark that (1) is equivalent to ${\bf j} = \{1,\ldots, p\}$ and 
      (2) is equivalent to ${\bf j} = \{q+1,\ldots, p+q\}$. Since these cases are products of a Hermitian symmetric domain and the cycle, by Proposition~\ref{pseudoconvex} $\ind(-\Lambda_{\ext}(E_0^{w_{\bf j}}))=0$.}\\

 Let $w_{\bf j} \in W_1^{\theta}$ be  as above.

$\bullet$ If $j_1 <j_r$, then $\lambda_+ = \ve_1 - \ve_r \in \Lambda_{\max}(E_0^{w_{\bf j}})$ and $$\{w_{\bf j}^{-1}(w\lambda_+): w \in W, w_{\bf j}^{-1}(w\lambda_+)>0\} =\{\ve_{j_a} - \ve_{j_b}:{a\leq p < b}, j_a <j_b\}.$$
 Let $(a_+, b_+)$ be such that $   b_+ -a_+=\min\{ b-a : {a\leq p < b}, j_a <j_b\}$. Then $w_{\max,+}\lambda_+ =\ve_{a_+} -\ve_{b_+}$ and  $\ind(-w_{\max,+}\lambda_+)=b_+ - a_+ -1$.
 {Remark that $b_+-a_+ -1 =0$ when ${\bf j} = \{1,\ldots, p\}$}.

$\bullet$ If $j_{p+1}<j_p$, then $\lambda_- = \ve_{p+1} - \ve_p \in \Lambda_{\max}(E_0^{w_{\bf j}})$ and
 $$\{w_{\bf j}^{-1}(w\lambda_-): w \in W, w_{\bf j}^{-1}(w\lambda_-)>0\} =\{\ve_{j_a} - \ve_{j_b}:b\leq p <p+1 \leq a, j_a <j_b\}.$$
 Let $(a_-, b_-)$ be such that $   a_- -b_-=\max\{ a-b : b\leq p<p+1 \leq a, j_a <j_b\}$. Then $w_{\max,-}\lambda_- =\ve_{a_-} -\ve_{b_-}$ and $\ind(-w_{\max,-}\lambda_-) = r-1-(a_--b_-)$.
{Remark that $r-1-(a_- -b_-) =0$ when ${\bf j} = \{q+1,\ldots, p+q\}$}.\\

{As a result, we have
 \begin{eqnarray*}
 \ind(-\Lambda_{\ext}(E_0^{w_{\bf j}}))=
 \min\{b_+ - a_+ -1, r-1-(a_- - b_-)  \}.
 \end{eqnarray*}}
 %Let $w_{\bf j} \in W_1^{\theta}$. Consider all pairs $(k_i, k_j)$ with $i<j$. Let $(i_+, j_+)$ and $(i_-, j_-)$ be such that
 %\begin{eqnarray*}
 %\vert k_{i_+} -k_{j_+}\vert &=&\min\{\vert k_i -k_j\vert :i<j, k_i <k_j\} \\
 %k_{i_-} -k_{j_-}&=&\max\{k_i-k_j: i<j, k_i>k_j \}.
 %\end{eqnarray*}
 %Then $w_{\max}\lambda_+ = \ve_{k_{i_+}} - \ve_{k_{j_+}}$ and $w_{\max}\lambda_-=\ve_{i_-} - \ve_{j_-}$ (For, when we act a simple reflection  in $W(\frak k, \frak t)$   on $\ve_i - \ve_j$, the difference $|i-j|$ is decreasing if $i<j$, and increasing if $i>j$).
 %$\ind(-w_{\max}\lambda_+)=|k_{i_+} -k_{j_+}|-1$ and %$\ind(-w_{\max}\lambda_-)=r-1-(k_{i_-}-k_{j_-})$.
 %
 %\begin{eqnarray*}
 %\ind(-\Lambda_{\ext}(E_0^{w_{\bf j}})) =\min\{ |k_{i_+} -k_{j_+}|-1, r-1-(k_{i_-}-k_{j_-})\}.
 %\end{eqnarray*}

%  For example,  if $p=2,q=5,r=7$ and the flag domain is  $D_{(+,-,+,-,+,+,+)}$, then $(k_1, \dots, k_r) = (3,1,4,2,5,6,7)$, and thus $(k_{i_+}, k_{j_+})=(2,5)$ and $(k_{i_-}, k_{j_-})=(3,1)$, and   $ \ind(-\Lambda_{\ext}(E_0^{w_{\bf j}}))=\min\{3-1, 7-1-2\}=2$.
\begin{example}
If $\frak g_0 = \frak{su}(3,4)$ and ${\bf j}=\{2,3,5\}$, then $(j_1, \dots, j_r) =(2,3,5,1,4,6,7)$ and $j_1<j_r$ and $j_{p+1}<j_p$. Since $(a_+, b_+)=(2,5)$ and $(a_-,b_-)=(4,1)$, we have $ \ind(-\Lambda_{\ext}(E_0^{w_{\bf j}}))= \min\{5-2-1, 7-1-(4-1)\} =2$.
\end{example}

 \subsection{$ \frak g_0=\frak{sp}(r, \mathbb R)$} \label{sp(r)}
{In this case, $G_0$ is of Hermitian type.}
Consider the skew-symmetric matrix
$$
J_{r,r} := \left(\begin{array}{ccccc}
&&&&-1\\
&&&1&\\
&&\cdot^{{\cdot}^{\cdot}}&&\\
&-1&&&\\
1&&&&\\
\end{array}\right)\in M_{2r}(\mathbb C).
$$
The Lie subalgebra $\frak g=\frak{sp}(r,\mathbb C)$ is the set of $X\in M_{2r}(\mathbb C)$ such that $J_{r,r}X+X^t J_{r,r}=0$ and $\frak b^{\text{ref}} = \frak{sp}(r,\mathbb C)\cap B^-_{2r}$.
{The associated flag manifold $G/B$ consists of all full flags
$$
0\subset F_1\subset F_2\subset\cdots\subset F_{r-1} \subset F_r 
$$
such that $F_j$ are isotropic subspaces in $\mathbb C^{2r}$ with respect to $J_{r,r}$ with 
$\dim_{\mathbb C} F_j = j$  for all $j=1,\ldots, r$.
Let $I_{r,r}$ be the standard Hermitian form of signature $(r,r)$ defined by the matrix $\left(\begin{array}{cc}
-I_r&0\\
0&I_r
\end{array}\right)$.
Let $a$ be a sequence $0\leq a_1\leq\cdots\leq a_r\leq r$ and $b$ be a sequence ${0}\leq b_1\leq \cdots \leq b_r\leq r$ such that $a_j+b_j=j$ for each $j=1,\ldots, r$. 
For a fixed pair $(a,b)$ the open orbit $D_{a,b}$ in $G/B$ is described by the set of all full flags $(0\subset F_1\subset\cdots\subset F_{r-1}\subset \mathbb C^r)$ such that $F_j$ has signature $(a_j,b_j)$ with respect to $I_{r,r}$.
There are $2^r$ number of flag domains in $G/B$ and each orbit is determined by the signature (cf. Chapter 5.5B in \cite{FHW}).
}

%The Lie algebra $\frak g_0=\frak{sp}(r,\mathbb R)=\frak{sp}(r,\mathbb C)\cap M_{2r}(\mathbb R)$ is a real form of $\frak g$. 
The subalgebra $\frak k$ is given by
$$
\left\{ \left( \begin{array}{cc}
A&0\\
0& J_{r,r}(-A^t) J_{r,r}
\end{array}\right) : A\in M_r(\mathbb C)\right\}
$$
and $\frak k \cap \frak{sp}(r,\mathbb R)\cong \frak u(r)$.
The Cartan subalgebras are given by $$\frak h = \frak t = \{\text{diag}(\ve_1,\ldots, \ve_r, -\ve_r,\ldots, -\ve_1): \ve_j\in \mathbb C \text{ for all } j\}.$$
The reference Borel subalgebra is given by
$\frak b_{\frak k}^{\text{ref}} = \frak{gl}(r,\mathbb C)\cap B^-_{2r}$.
{The longest root of $\Sigma(\frak g, \frak h)$ is $\lambda_+$. }

 Positive noncompact roots of $(\frak h, \frak g)$ are {$\{\varepsilon_i+\varepsilon_j\}_{1\leq i
<j\leq r}\cup \{2\varepsilon_i\}_{1\leq i\leq r}$}
and positive compact roots of $(\frak h, \frak g)$ are {\{$\varepsilon_i-\varepsilon_j\}_{1\leq i<j\leq r}$}.
The Hasse diagram of the orbit $W.\lambda_+=\{w\lambda_+:w \in W\}$ is
 \begin{center}
 \begin{tikzcd}[column sep=0.7cm, nodes={anchor=center}]
\lambda_+=2\ve_1 \rar& 2\ve_2\rar&\cdots \rar& 2 \ve_{r-2}\rar & 2 \ve_{r-1} \rar& 2 \ve_r=-\lambda_-
  \end{tikzcd} 
\end{center}
 and
that of $W.\lambda_-=\{w\lambda_-:w \in W\}$ is
 \begin{center}
 \begin{tikzcd}[column sep=0.7cm, nodes={anchor=center}]
\lambda_-=-2\ve_r \rar& -2 \ve_{r-1}\rar & \cdots \rar& -2 \ve_{2} \rar& -2 \ve_1=-\lambda_+ .
  \end{tikzcd} 
\end{center}
Since we have $\lambda_+=2 \ve_1$ and $w_{\max}\lambda_+ = 2 \ve_r = w_0 \lambda_+$, we obtain that $-w_{\max}\lambda_+$ is dominant and hence $\ind (-w_{\max}\lambda_+)=0$. \\

In this case, we have $W_1^{\theta} =\left\{ w_{\bf j}: {\bf j} \subset \{1, \dots, r\} \right\}$ and thus there are $2^r$ base cycles in $G/B$.  For ${\bf j}=\{j_1, \dots, j_p\} \subset \{1, \dots, r\}$ with $j_1 <j_2 < \dots <j_p$, after ordering the complement $\{1,2, \dots, r\} \backslash {\bf j}$ in increasing order, as $j_{p+1}<\dots<j_{r}$,  define $w_{\bf j} \in W^{\theta}$ by
 \begin{eqnarray*}
 w_{\bf j}^{-1}(\ve_{ i}) =\left\{\begin{array}{ll}
 \ve_{j_i} & \text{ for } 1 \leq i \leq p \\
 -\ve_{j_{r+p+1-i}} & \text{ for } p+1 \leq i \leq r
 \end{array}\right.
 \end{eqnarray*}
 for $1 \leq i \leq r$.
 For $w_{\bf j} \in W_1^{\theta}$,
 $w_{\bf j}^{-1}(w\lambda_+), w \in W$ is
 \begin{center}
 \begin{tikzcd}[column sep=0.5cm, nodes={anchor=center}]
 2\ve_{j_1} &2\ve_{j_2}& \cdots & 2 \ve_{j_{p}} &-2\ve_{j_{r}} & \cdots & -2 \ve_{j_{p+1}} 
  \end{tikzcd} 
\end{center}
 and $w_{\bf j}^{-1}(w\lambda_-), w \in W$ is
 \begin{center}
 \begin{tikzcd}[column sep=0.5cm, nodes={anchor=center}]
 2\ve_{j_{p+1}} & \dots &  2 \ve_{j_{r}} &-2\ve_{j_{p}} & \dots & -2 \ve_{j_{ 1}} .
  \end{tikzcd} 
\end{center}

 %For $w_{\bf j} \in W_1^{\theta}$ with $|{\bf j}|=p$, write  $w_{\bf j}(\ve_i) =\delta_i\ve_{k_i}$ for $1 \leq i \leq r$, where $\delta_i \in \{1,-1\}$. Then $\delta_i=1$ if $i \in {\bf j}$ and $\delta_i=-1$ if $i \not\in {\bf j}$. Furthermore,
 %$$\Lambda( \frak s \cap w_{\bf j}(\frak b^{{\rm ref},n})) \cap W.\{\lambda_+,\lambda_-\}=\{2\delta_1 \ve_{k_1}, \dots, 2\delta_p\ve_{k_p}, 2\delta_{p+1}\ve_{k_{p+1}}, \dots,  2\delta_{r}\ve_{k_r}\}$$
 % and thus
 % $$ \Lambda( \frak s \cap w_{\bf j}(\frak b^{{\rm ref},n} )) \cap W.\{\lambda_+,\lambda_-\}=\{2\ve_1, \dots, 2\ve_p, -2\ve_{p+1}, \dots, -2\ve_r\}.$$
\noindent Thus $w_{\max,+}\lambda_+ = 2 \ve_p$ and $\ind(-w_{\max,+}\lambda_+) =r-p$ and $w_{\max,-}\lambda_- = -2\ve_{p+1}$ and $\ind(-w_{\max,-}\lambda_-)=p$. Hence
 $\ind(-\Lambda_{\ext}(E_0^{w_{\bf j}}))= \min\{p, r-p\}$.

%\medskip
 \subsection{$\frak g_0=\frak{so}(m, n)$}\label{so(m,n)}
 In this case only when $m=2$, $G_0$ is of Hermitian type. Let $S$ be a nondegenerate symmetric bilinear form on $\mathbb C^{m+n}$, which is defined by the antidiagonal matrix with antidiagonal entries to be $1$. Then $SO(m+n, \mathbb C)$ is a set of elements in $SL(m+n,\mathbb C)$ preserving $S$.
We have $\frak g = \frak{so}(m+n,\mathbb C)$ and $\frak b^{\text{ref}}=\frak{so}(m+n,\mathbb C)\cap B_{m+n}^-$. Let $I_{m,n}$ be the standard Hermitian form of signature $(m,n)$ defined by the matrix $
\left(\begin{array}{cc}
-I_n&0\\
0&I_m
\end{array}\right)$.
{Let $p=\left[\frac{m}{2}\right]$, $q=\left[\frac{n}{2}\right]$ and $r=p+q$.
The associated flag manifold $G/B$ consists of all full flags
$$
0\subset F_1\subset F_2\subset\cdots\subset F_{r-1} \subset F_r \subset  \mathbb C^{m+n}
$$
such that $F_j$ are isotropic subspace with respect to $S$ of $\dim F_j = j$ for all $j=1,\ldots, r$.
Let $a$ be a sequence $0\leq a_1\leq\cdots\leq a_r\leq n$ and $b$ be a sequence ${0}\leq b_1\leq \cdots \leq b_r\leq m$ such that $a_j+b_j=j$ for each $j=1,\ldots, n+m$.} 
For a fixed pair $(a,b)$ the open orbit $D_{a,b}$ in $G/B$ is described by the set of all full flags $(0\subset F_1\subset\cdots\subset F_{r-1}\subset F_{r} \subset \mathbb C^{m+n})$ such that $F_j$ has signature $(a_j,b_j)$ with respect to $I_{m,n}$. 
There are $\frac{(p+q)!}{p!q!}$ number of flag domains in $G/B$ and each orbit is determined by the signature.

We have $\text{rank}(\frak g) = \text{rank}(\frak k)$ except for the case when $m$ and $n$ are both odd. Moreover, $\lambda_{\frak s}$ is not induced from the longest root of $\Sigma(\frak g, \frak h)$ except for the case when $G_0$ is of Hermitian type (i.e., $m=2$).
\medskip

 \noindent \textsf{\textbf{Case 1.}} If $(m,n)=(2p+1, 2q+1)$ with ${0}\leq p \leq q$, then $\frak g = \frak{so}(2(r+1), \mathbb C)$ with $r=p+q$
and $\frak k = \frak{so}(2p+1,\mathbb C)\times \frak{so}(2q+1,\mathbb C)$.
Consider $\frak g$ as a Lie subalgebra in $\frak{sl}(2(r+1),\mathbb C)$ as in Section~\ref{sl}. Then
$\frak h = \text{diag}(\ve_1, \ldots,
\ve_{r+1}, -\ve_{r+1},\ldots, -\ve_1)$ and $\frak t = \text{diag}(\ve_1, \ldots,
\ve_r,0,0, -\ve_{r},\ldots, -\ve_1)$.
The reference Borel subalgebras $\frak b^{\text{ref}}$ and $\frak b_{\frak k}^{\text{ref}}$ are given by $\frak{so}(2(r+1),\mathbb C)\cap B_{2(r+1)}^-$ and $\left(\frak{so}(2p+1,\mathbb C)\times \frak{so}(2q+1,\mathbb C)\right) \cap B^-_{2(r+1)}$, respectively.
The simple root system of $\frak k$ is given by $\Psi_{\frak k}^{(1)} \cup \Psi_{\frak k}^{(2)}$ where
\begin{equation}\nonumber
\begin{aligned}
\Psi_{\frak k}^{(1)} &= \{ \ve_1-\ve_2, \ldots, \ve_{p-1}-\ve_p, \ve_p\}=: \{ \beta_1^{(1)},\ldots, \beta_{p}^{(1)}\}\quad\text{ and } \\
\Psi_{\frak k}^{(2)} &= \{ \ve_{p+1}-\ve_{p+2}, \ldots, \ve_{r-1}-\ve_r, \ve_r\}=: \{ \beta_1^{(2)},\ldots, \beta_{q}^{(2)}\}.
\end{aligned}
\end{equation}
{Remark that if $p=0$, $\Psi_{\frak k}^{(1)}$ is an empty set and $\Psi_{\frak k}^{(2)}=\{\ve_1-\ve_2,\ldots, \ve_{r-1}-\ve_r, \ve_r\}$.}
{For $p\geq 1$,} since $\lambda_{\frak s} =\xi_1^{(1)} + \xi_1^{(2)} = \ve_1 + \ve_{p+1}$, we have $$W.\lambda_{\frak s} =\{ \ve_a + \ve_b, \ve_a - \ve_b, -\ve_a+\ve_b, -\ve_a -\ve_b: a   \leq p <p+1 \leq b  \}.$$

\begin{figure}[h]
  {  \footnotesize\begin{equation} \nonumber  
\arraycolsep=1.4pt\def\arraystretch{2}
\begin{array}{|cccc|cccc|}
\hline
\lambda_{\frak s} = \ve_{1}+\ve_{p+1}& \ve_{1}+\ve_{p+2} &\cdots &\ve_{1}+\ve_{p+q}&\ve_{1}-\ve_{p+q}&\ve_{1}-\ve_{p+q-1}&\cdots&\ve_{1}-\ve_{p+1}\\
 \ve_{2}+\ve_{p+1}& \ve_{2}+\ve_{p+2} &\cdots &\ve_{2}+\ve_{p+q}&\ve_{2}-\ve_{p+q}&\ve_{2}-\ve_{p+q-1}&\cdots&\ve_{2}-\ve_{p+1}\\
 \vdots& \vdots& \vdots& \vdots& \vdots& \vdots& \vdots& \vdots\\
  \ve_{p}+\ve_{p+1}& \ve_{p}+\ve_{p+2} &\cdots &\ve_{p}+\ve_{p+q}&\ve_{p}-\ve_{p+q}&\ve_{p}-\ve_{p+q-1}&\cdots&\ve_{p}-\ve_{p+1}\\\hline
- \ve_{p}+\ve_{p+1}& -\ve_{p}+\ve_{p+2} &\cdots &-\ve_{p}+\ve_{p+q}&-\ve_{p}-\ve_{p+q}&-\ve_{p}-\ve_{p+q-1}&\cdots&-\ve_{p}-\ve_{p+1}\\
 -\ve_{p-1}+\ve_{p+1}&  -\ve_{p-1}+\ve_{p+2} &\cdots & -\ve_{p-1}+\ve_{p+q}& -\ve_{p-1}-\ve_{p+q}& -\ve_{p-1}-\ve_{p+q-1}&\cdots& -\ve_{p-1}-\ve_{p+1}\\
 \vdots& \vdots& \vdots& \vdots& \vdots& \vdots& \vdots& \vdots\\
  -\ve_{1}+\ve_{p+1}& -\ve_{1}+\ve_{p+2} &\cdots &-\ve_{1}+\ve_{p+q}&-\ve_{1}-\ve_{p+q}&-\ve_{1}-\ve_{p+q-1}&\cdots&-\ve_{1}-\ve_{p+1}\\\hline
\end{array}
\end{equation}}
\caption{{Hasse diagram for} $\frak{so}(2p+1,2q+1)$}
\label{oddodd}
\end{figure}

\noindent{Figure~\ref{oddodd} is a Hasse diagram similar to Figure~\ref{su(p,q)}, but with the arrows omitted. For each $\ve_i-\ve_j$, arrows are directed to the downward and rightward. The roots in the first and second quadrants are positive, while those in the third and fourth quadrants are negative.
}

Since for ${a\leq p < b}$,
\begin{equation}\nonumber
\begin{aligned}
 \ind(\ve_a -\ve_b)&=(a-1)+(p+q-b) +q, \text{ and } \\
 \ind(-\ve_a +\ve_b) &= (b-(p+1)) + (p-a) +p ,
\end{aligned}
\end{equation}
we obtain
 $w_{\max}\lambda_{\frak s} = \ve_p -\ve_{p+1}$. Thus $-w_{\max}\lambda_{\frak s}=-\ve_p +\ve_{p+1}$ and
hence $\ind(-w_{\max}\lambda_{\frak s}) = p$. \\

{If $p=0$, since $\lambda_{\frak s} = \xi_1= \ve_1$, $W.\lambda_{\frak s}$ is
 \begin{eqnarray*}
  \begin{array}{cccccccccc}
  \ve_1& \ve_2& \cdots& \ve_{r-1}& \ve_r& -\ve_r& -\ve_{r-1}& \cdots& -\ve_2& -\ve_1
  \end{array}
  \end{eqnarray*}
  and $w_{\max}\lambda_{\frak s} = \ve_r$. Thus $-w_{\max}\lambda_{\frak s}=-\ve_r$ and
hence $\ind(-w_{\max}\lambda_{\frak s}) = r=q$.}

\vspace{0.2cm}
In this case, we have $W_1^{\theta} =\{ w_{\bf j}: {\bf j} \subset \{1, \dots, r\}, |{\bf j}|=p\}$, and thus there exists $\frac{r!}{p!q!}$ base cycles in $G/B$.
For a subset ${\bf j}=\{j_1, \dots, j_p\} \subset \{1, \dots, r\}$ with $j_1 <j_2 < \dots <j_p$, after ordering the complement $\{1,2, \dots, r\} \backslash {\bf j}$ in increasing order, as $j_{p+1}<\dots<j_{r}$,  define $w_{\bf j} \in W^{\theta}$ by
 $$w_{\bf j}^{-1}(\ve_{ i}) =\ve_{j_i}$$
 for $i = 1,\ldots, r$.
 For any $\ve_a + \ve_b \in W.\lambda_{\frak s}$, we have $w_{\bf j}^{-1}(\ve_a+ \ve_b)=\ve_{j_a} + \ve_{j_b} \in \Lambda(\frak b^{{\rm ref},n})$
  and for any $-\ve_a - \ve_b \in W.\lambda_{\frak s}$, we have $w_{\bf j}^{-1}(-\ve_a- \ve_b)=-\ve_{j_a} - \ve_{j_b} \not\in \Lambda(\frak b^{{\rm ref},n})$.
 It suffices to check for which $\ve_a -\ve_b \in W.\lambda_{\frak s}$, $w_{\bf j}^{-1}(\ve_a- \ve_b)=\ve_{j_a} -\ve_{j_b} \in \Lambda(\frak b^{{\rm ref},n})$.

\medskip

 Let $(a', b')$ and $(a'', b'')$ be such that
 {
\begin{equation}\label{ab}
\begin{aligned}
 b'-a'&=\left\{
\begin{array}{cc}
\min\{ b-a: a \leq p <p+1 \leq b, j_a <j_b \}&\text{ if }\, {\bf j} \neq \{q+1, \ldots, p+q\},\\
0&\text{ if } \,{\bf j} = \{q+1, \ldots, p+q\},\\
\end{array}\right.\\
 b'' - a''&=\left\{
\begin{array}{cc}
     { \max}\{ b-a: a \leq p <p+1 \leq b, j_a >j_b\} 
     &\text{ if }\, {\bf j} \neq \{1,\ldots, p\},\\
     0 
     &\text{ if }\, {\bf j} = \{1,\ldots, p\}.\\
\end{array}\right.
\end{aligned}
\end{equation}}
 Then
 $\ind(-\Lambda_{\ext}(E_0^{w_{\bf j}}))= \min\{(b'-(p+1)) + (p-a') +p, (a''-1)+ (p+q-b'') +q \}$.

\begin{example} $\,$

\begin{enumerate}
\item[$\bullet$] When $\frak g_0 = \frak{so}(7,9)$, i.e.
$(p,q) = (3,4)$ and ${\bf j }= \{2,5,6\}$,
 $w_{\bf j}^{-1}(W.\lambda_{\frak s})$ is given by Figure~\ref{so(7,9)_1}
 and $w_{\text{max}}\lambda_{\frak s} = \ve_1-\ve_5, \ve_3-\ve_7$, and $-\ve_2+\ve_6$.
  Hence $\ind(-w_{\max}\lambda_{\frak s}) = 6$.

\item [$\bullet$]
When ${\bf j }= \{4,6,7\}$,
 $w_{\bf j}^{-1}(W.\lambda_{\frak s})$ is given by Figure~\ref{so(7,9)_2}
  and {$w_{\text{max}}\lambda_{\frak s} = -\ve_1+\ve_6, -\ve_2+\ve_7$.}
  Hence $\ind(-w_{\max}\lambda_{\frak s}) = 5$.

\end{enumerate}
\end{example}

Here we explain notations in   Figure~\ref{so(7,9)_0} and  Figure~\ref{so(7,9)_1} and  Figure~\ref{so(7,9)_2}. For simplicity, we write $\ve_i-\ve_j$, $\ve_i+\ve_j$, $-\ve_i+\ve_j$ and $-\ve_i-\ve_j$ as $(i,-j)$, $(i,j)$, $(-i, j)$ and $(-i,-j)$, respectively. Each root in Figure~\ref{so(7,9)_0} is mapped by $w_{\bf j}^{-1}$ to the element in the same entry in Figure~\ref{so(7,9)_1}.
The roots highlighted in blue color in Figure~\ref{so(7,9)_1} represent the positive roots with respect to $\frak b^{\text{ref},n}$, and the roots highlighted in blue color in Figure~\ref{so(7,9)_0} represent the positive roots  with respect to  $w_{\bf j}(\frak b^{\text{ref},n})$.

   \begin{figure}[h]
 { \footnotesize \begin{equation}\nonumber
\arraycolsep=2pt\def\arraystretch{2} \begin{array}{|cccc|cccc|}\hline
{\color{blue}(1,4)} & {\color{blue}(1,5)}& {\color{blue}(1,6)}& {\color{blue}(1,7)}& {\color{blue}(1,-7)}&  {\color{blue}(1,-6)}& {\color{blue}(1,-5)}& (1,-4)\\
{\color{blue}(2,4)}& {\color{blue}(2,5)}& {\color{blue}(2,6)}& {\color{blue} (2,7)}&  {\color{blue}(2,-7)}& (2,-6)&(2,-5)& (2,-4)\\
{\color{blue}(3,4)}& {\color{blue}(3,5)}& {\color{blue}(3,6)}& {\color{blue}(3,7)}&  {\color{blue}(3,-7)}& (3,-6)&(3,-5)& (3,-4)\\\hline
 {\color{blue}(-3,4)}& {\color{blue} (-3,5)}&  {\color{blue}(-3,6)}& (-3,7)& (-3,-7)& (-3,-6)&(-3,-5)& (-3,-4)\\
 {\color{blue}(-2,4)}&  {\color{blue}(-2,5)}& {\color{blue} (-2,6)}& (-2,7)& (-2,-7)& (-2,-6)&(-2,-5)& (-2,-4)\\
 {\color{blue}(-1,4)}& (-1,5)& (-1,6)& (-1,7)& (-1,-7)& (-1,-6)&(-1,-5)& (-2,-4)\\\hline
 \end{array}
 \end{equation}}
 \caption{{Hasse diagram for} $\frak{so}(7,9)$} 
 \label{so(7,9)_0}
 \end{figure}

 \begin{figure}[h]
 { \footnotesize \begin{equation}\nonumber
\arraycolsep=2pt\def\arraystretch{2} \begin{array}{|cccc|cccc|}\hline
{\color{blue}(2,1)} & {\color{blue}(2,3)}& {\color{blue}(2,4)}& {\color{blue}(2,7)}& {\color{blue}(2,-7)}&  {\color{blue}(2,-4)}& {\color{blue}(2,-3)}& (2,-1)\\
{\color{blue}(5,1)}& {\color{blue}(5,3)}& {\color{blue}(5,4)}& {\color{blue} (5,7)}&  {\color{blue}(5,-7)}& (5,-4)&(5,-3)& (5,-1)\\
{\color{blue}(6,1)}& {\color{blue}(6,3)}& {\color{blue}(6,4)}& {\color{blue}(6,7)}&  {\color{blue}(6,-7)}& (6,-4)&(6,-3)& (6,-1)\\\hline
 {\color{blue}(-6,1)}& {\color{blue} (-6,3)}&  {\color{blue}(-6,4)}& (-6,7)& (-6,-7)& (-6,-4)&(-6,-3)& (-6,-1)\\
 {\color{blue}(-5,1)}&  {\color{blue}(-5,3)}& {\color{blue} (-5,4)}& (-5,7)& (-5,-7)& (-5,-4)&(-5,-3)& (-5,-1)\\
 {\color{blue}(-2,1)}& (-2,3)& (-2,4)& (-2,7)& (-2,-7)& (-2,-4)&(-2,-3)& (-2,-1)\\\hline
 \end{array}
 \end{equation}}
 \caption{Hasse diagram for} {$\frak{so}(7,9)$ with ${\bf j} = \{2,5,6\}$}
 \label{so(7,9)_1}
 \end{figure}

 \begin{figure}[h]
 {  \footnotesize \begin{equation}\nonumber
\arraycolsep=2pt\def\arraystretch{2} \begin{array}{|cccc|cccc|}\hline
{\color{blue}(4,1)}&{\color{blue}(4,2)}&{\color{blue}(4,3)}&{\color{blue}(4,5)}&{\color{blue}(4,-5)}&{(4,-3)}&{(4,-2)}&{(4,-1)}\\
{\color{blue}(6,1)}&{\color{blue}(6,2)}&{\color{blue}(6,3)}&{\color{blue}(6,5)}&{(6,-5)}&{(6,-3)}&{(6,-2)}&{(6,-1)}\\
{\color{blue}(7,1)}&{\color{blue}(7,2)}&{\color{blue}(7,3)}&{\color{blue}(7,5)}&{(7,-5)}&{(7,-3)}&{(7,-2)}&{(7,-1)}\\\hline
{\color{blue}(-7,1)}&{\color{blue}(-7,2)}&{\color{blue}(-7,3)}&{\color{blue}(-7,5)}&{(-7,-5)}&{(-7,-3)}&{(-7,-2)}&{(-7,-1)}\\
{\color{blue}(-6,1)}&{\color{blue}(-6,2)}&{\color{blue}(-6,3)}&{\color{blue}(-6,5)}&{(-6,-5)}&{(-6,-3)}&{(-6,-2)}&{(-6,-1)}\\
{\color{blue}(-4,1)}&{\color{blue}(-4,2)}&{\color{blue}(-4,3)}&{(-4,5)}&{(-4,-5)}&{(-4,-3)}&{(-4,-2)}&{(-4,-1)}\\\hline
 \end{array}
 \end{equation}}
\caption{{Hasse diagram for} $\frak{so}(7,9)$ with ${\bf j} = \{4,6,7\}$}
\label{so(7,9)_2}
\end{figure}

 \noindent \textsf{\textbf{Case 2.}} If $(m,n)=(2p, 2q+1)$ with $p > 2$, $q\geq {0}$, then $\frak g = \frak{so}(2r+1, \mathbb C)$ with $r=p+q$
and $\frak k = \frak{so}(2p,\mathbb C)\times \frak{so}(2q+1,\mathbb C)$.
In this case
$\frak h =\frak t = \{ \text{diag}(\ve_1, \cdots,
\ve_r,0, -\ve_{r},\ldots, -\ve_1): \ve_j\in \mathbb C \text{ for all } j\}$.
The reference Borel subalgebras $\frak b^{\text{ref}}$ and $\frak b_{\frak k}^{\text{ref}}$ are given by $\frak{so}(2r+1,\mathbb C)\cap B_{2(r+1)}^-$ and $\left(\frak{so}(2p,\mathbb C)\times \frak{so}(2q+1,\mathbb C)\right) \cap B^-_{2(r+1)}$, respectively.
The simple root system of $\frak k$ is given by $\Psi_{\frak k}^{(1)} \cup \Psi_{\frak k}^{(2)}$ where
\begin{equation}\nonumber
\begin{aligned}
\Psi_{\frak k}^{(1)} &= \left\{ \ve_1-\ve_2, \ldots, \begin{array}{c}
\ve_{p-1}-\ve_p\\
\ve_{p-1}+ \ve_p
\end{array}
\right\}
=: \left\{ \beta_1^{(1)},\ldots, \begin{array}{c}
\beta_{p-1}^{(1)}\\
\beta_p^{(1)}
\end{array}
\right\}\quad\text{ and } \\
\Psi_{\frak k}^{(2)} &= \{ \ve_{p+1}-\ve_{p+2}, \ldots, \ve_{r-1}-\ve_r, \ve_r\}=: \{ \beta_1^{(2)},\ldots, \beta_{q}^{(2)}\}.
\end{aligned}
\end{equation}
If $q\geq 1$, since $\lambda_{\frak s} =\xi_1^{(1)} + \xi_1^{(2)} = \ve_1 + \ve_{p+1}$, we have
$$W. \lambda_{\frak s} =\{\ve_a +\ve_b, \ve_a -\ve_b, {-\ve_a+\ve_b}, -\ve_a -\ve_b: a \leq p<p+1 \leq b\},$$ the set itself is the same as the case of $(m,n) = (2p+1, 2q+1)$. 
%\noindent
However the indices are different; for ${a\leq p < b}$,
 \begin{eqnarray*}
 \ind(\ve_a -\ve_b)&=& (a-1)+(p+q-b) +q , \\
 \ind(-\ve_a +\ve_b) &=& (b-(p+1)) + (p-a) +p -1 
 \end{eqnarray*}
 {(See Figure~\ref{evenodd} and Figure~\ref{partial Hasse so(2p,2q+1)}).}
 Thus {$w_{\textup{max}} \lambda_{\frak s} = \ve_p-\ve_{p+1}$ and hence} $\ind(-\Lambda_{\ext}(E_0))=p-1$.

%\pgfimage[width=10cm]{diagram 2p2q+1}
%\includegraphics[width=\textwidth]{figures/your-figure.eps}
%\includegraphics[width=\textwidth]{figures/diagram 2p2q+1.png}

%********************************** diagram 2p2q+1
%\begin{figure}[!h]

 %       \centering

 %       \includegraphics[angle=0, width=1.5\textwidth]{diagram 2p2q+1.png}

 %       \caption{so(2p,2q+1)}

 %       \label{fig1}

%\end{figure}
%*********************************************

\begin{figure}[h]
{   \footnotesize\begin{equation}\nonumber
\arraycolsep=1.4pt\def\arraystretch{2}\begin{array}{|cccc|cccc|}
\hline
\lambda_{\frak s} = \ve_{1}+\ve_{p+1}& \ve_{1}+\ve_{p+2} &\cdots &\ve_{1}+\ve_{p+q}&\ve_{1}-\ve_{p+q}&\ve_{1}-\ve_{p+q-1}&\cdots&\ve_{1}-\ve_{p+1}\\
 \ve_{2}+\ve_{p+1}& \ve_{2}+\ve_{p+2} &\cdots &\ve_{2}+\ve_{p+q}&\ve_{2}-\ve_{p+q}&\ve_{2}-\ve_{p+q-1}&\cdots&\ve_{2}-\ve_{p+1}\\
 \vdots& \vdots& \vdots& \vdots& \vdots& \vdots& \vdots& \vdots\\
  \ve_{p-1}+\ve_{p+1}& \ve_{p-1}+\ve_{p+2} &\cdots &\ve_{p-1}+\ve_{p+q}&\ve_{p-1}-\ve_{p+q}&\ve_{p-1}-\ve_{p+q-1}&\cdots&\ve_{p-1}-\ve_{p+1}\\\hline
\pm \ve_{p}+\ve_{p+1}& \pm\ve_{p}+\ve_{p+2} &\cdots &\pm\ve_{p}+\ve_{p+q}& \pm\ve_{p}-\ve_{p+q}&\pm\ve_{p}-\ve_{p+q-1}&\cdots&\pm\ve_{p}-\ve_{p+1}\\\hline
 -\ve_{p-1}+\ve_{p+1}&  -\ve_{p-1}+\ve_{p+2} &\cdots & -\ve_{p-1}+\ve_{p+q}& -\ve_{p-1}-\ve_{p+q}& -\ve_{p-1}-\ve_{p+q-1}&\cdots& -\ve_{p-1}-\ve_{p+1}\\
 \vdots& \vdots& \vdots& \vdots& \vdots& \vdots& \vdots& \vdots\\
  -\ve_{1}+\ve_{p+1}& -\ve_{1}+\ve_{p+2} &\cdots &-\ve_{1}+\ve_{p+q}&-\ve_{1}-\ve_{p+q}&-\ve_{1}-\ve_{p+q-1}&\cdots&-\ve_{1}-\ve_{p+1}\\\hline
\end{array}
\end{equation}}
\caption{{Hasse diagram for} $\frak{so}(2p, 2q+1)$}
\label{evenodd}
\end{figure}
\noindent{Figure~\ref{evenodd} is a Hasse diagram similar to Figure~\ref{su(p,q)}, but with the arrows omitted. For each $\ve_i\pm\ve_j$, $-\ve_i\pm\ve_j$ with $i\neq p-1, p$ arrows are directed to the downward and rightward. 
The middle row of Figure~\ref{evenodd}, i.e. $i=p-1, p$, represents the Hasse diagram Figure~\ref{partial Hasse so(2p,2q+1)}.
%The roots in the first and second quadrants are positive, while those in the third and fourth quadrants are negative.

\begin{figure}[h]
\begin{tikzpicture}[scale=0.85]
  \node (a) at (0,4) {$\ve_{p-1}+\ve_{p+1}$};
  \node (b) at (-1,3) {$\ve_{p}+\ve_{p+1}$};
  \node (c) at (1,2) {$-\ve_{p}+\ve_{p+1}$};
  \node (d) at (0,1) {$-\ve_{p-1}+\ve_{p+1}$};
  \node (e) at (4,4) {$\ve_{p-1}+\ve_{p+2}$};
  \node (f) at (3,3) {$\ve_{p}+\ve_{p+1}$};
  \node (g) at (5,2) {$-\ve_{p}+\ve_{p+1}$};
  \node (h) at (4,1) {$-\ve_{p-1}+\ve_{p+2}$};

  \node (i) at (7,4) {$\cdots$};
  \node (j) at (7,3) {$\cdots$};
  \node (k) at (7,2) {$\cdots$};
  \node (l) at (7,1) {$\cdots$};

  \node (m) at (10,4) {$\ve_{p-1}-\ve_{p+2}$};
  \node (n) at (9,3) {$\ve_{p}-\ve_{p+2}$};
  \node (o) at (11,2) {$-\ve_{p}-\ve_{p+2}$};
  \node (p) at (10,1) {$-\ve_{p-1}-\ve_{p+2}$};
  \node (q) at (14,4) {$\ve_{p-1}-\ve_{p+1}$};
  \node (r) at (13,3) {$\ve_{p}-\ve_{p+1}$};
  \node (s) at (15,2) {$-\ve_{p}-\ve_{p+1}$};
  \node (t) at (14,1) {$-\ve_{p-1}-\ve_{p+1}$};

  \draw[->](a)--(e);
  \draw[->](a)--(b);
  \draw[->](a)--(c);
  \draw[->](a)--(e);
  
  \draw[->](b)--(d);
  \draw[->](b)--(f);
  
  \draw[->](c)--(d);
  \draw[->](c)--(g);
  
  \draw[->](d)--(h);
  
  \draw[->](e)--(f);
  \draw[->](e)--(g);
  \draw[->](e)--(i);
    
  \draw[->](f)--(h);
  \draw[->](f)--(j);

  \draw[->](g)--(h);
  \draw[->](g)--(k);

  \draw[->](h)--(l);
  
  \draw[->](m)--(q);
  \draw[->](m)--(n);
  \draw[->](m)--(o);
  \draw[->](m)--(q);
  
  \draw[->](n)--(p);
  \draw[->](n)--(r);
  
  \draw[->](o)--(p);
  \draw[->](o)--(s);
  
  \draw[->](p)--(t);
  
  \draw[->](q)--(r);
  \draw[->](q)--(s);

  \draw[->](r)--(t);

  \draw[->](s)--(t);
  
  \draw[->](i)--(m);
  \draw[->](j)--(n);
  \draw[->](k)--(o);
  \draw[->](l)--(p);

\end{tikzpicture}
\caption{Partial Hasse diagram for $\frak{so}(2p,2q+1)$}
\label{partial Hasse so(2p,2q+1)}
\end{figure}
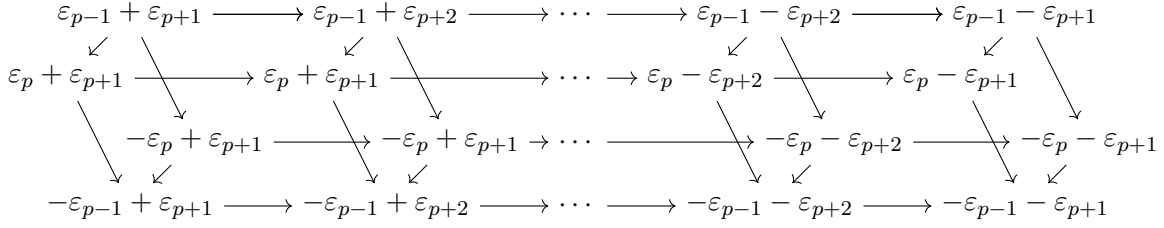
}

{If $q=0$, then $\Psi_{\frak k}^{(2)}$ is an empty set. Since $\lambda_{\frak s} = \xi_1=\ve_1$ and the Hasse diagram of $W.\lambda_{\frak s}$ is
\begin{center}
\begin{tikzpicture}
    \matrix (a) [matrix of math nodes, column sep=0.35cm ]{
    &&&&2\ve_r&&&\\
    \ve_1
     & \ve_2
     & \cdots
     & \ve_{r-1}
     &~
     & -\ve_{r-1}
     & \cdots
     & -\ve_2
     & -\ve_1.\\
     &&&&-\ve_r&&&\\
    };
    \draw[->] (a-2-1)--(a-2-2);
    \draw[->] (a-2-2)--(a-2-3);
    \draw[->] (a-2-3)--(a-2-4);
    \draw[->] (a-2-4)--(a-1-5);
    \draw[->] (a-2-4)--(a-3-5);
    \draw[->] (a-1-5)--(a-2-6);
    \draw[->] (a-3-5)--(a-2-6);
    \draw[->] (a-2-6)--(a-2-7);
    \draw[->] (a-2-7)--(a-2-8);
    \draw[->] (a-2-8)--(a-2-9);
\end{tikzpicture}
 \end{center}
Thus $w_{\max}\lambda_{\frak s} =  \ve_r$ and $\ind(-w_{\max}\lambda_{\frak s}) =r-1 =p-1$. \\
}
\medskip

Since we have $W_1^{\theta}
=\left\{ w_{\bf j}, w'_{\bf j}:{\bf j} \subset \{1, \dots, r\} \right\}$, where
  \begin{eqnarray*}
  w_{\bf j}^{-1}(\ve_i) &=&\ve_{j_i} \text{ for all } i, \text{ and }\\
  {w_{\bf j}'}^{-1}(\ve_i) &=& \left\{
  \begin{array}{cc}
  \ve_{j_i} & \text{ for } i \not=p,\\
  -\ve_{j_i} & \text{ for }i=p.
  \end{array} \right.
  \end{eqnarray*}
For ${\bf j} \subset \{1, \dots, r\}$, define $(a',b')$ and $(a'',b'')$ as in the case when $(m,n)=(2p+1, 2q+1)$. Then
$\ind(-\Lambda_{\ext}(E_0^{w_{\bf j}}))=\min\{(b'-(p+1)) + (p-a') +p-1, (a'' -1) + (p+q-b'')+q\}$.

  \medskip

 \noindent \textsf{\textbf{Case 3.}}  If $(m,n) =(2p, 2q)$ with $2\leq p\leq q$, {then $\frak g = \frak{so}(2r,\mathbb C)$ with $r=p+q$ and $\frak k\cong \frak{so}(2p)\times \frak{so}(2q)$. In this case, 
 $\frak h = \frak t = \{\textup{diag}(\ve_1, \ldots, \ve_r, -\ve_r,\ldots,-\ve_1) : \ve_j\in \mathbb C \textup{ for all } j\}$.
}
 We compute in a similar way to the case $(m,n)=(2p+1, 2q+1)$ and $(m,n)=(2p, 2q+1)$.
As a result we obtain
$$
W. \lambda_{\frak s} =\{\ve_a +\ve_b, \ve_a -\ve_b, -\ve_a + \ve_b, -\ve_a -\ve_b: a \leq p<p+1 \leq b\}.
$$
{The Hasse diagram corresponding to $W. \lambda_{\frak s}$ is shown in Figure~\ref{eveneven}. We construct it using a similar approach as in Figure~\ref{evenodd}.}

%\pgfimage[width=10cm]{diagram 2p2q+1}
%\includegraphics[width=\textwidth]{figures/your-figure.eps}
%\includegraphics[width=\textwidth]{figures/diagram 2p2q+1.png}

%****************************** diagram 2p2q
%\begin{figure}[!h]

 %       \centering

 %       \includegraphics[angle=0, width=1.5\textwidth]{diagram 2p2q.png}

%        \caption{so(2p,2q)}

 %       \label{fig2}

%\end{figure}
%**********************************
\begin{figure}[h]
{ \footnotesize \begin{equation} \nonumber
\arraycolsep=1.4pt\def\arraystretch{1.5}\begin{array}{|cccc|c|ccc|}
\hline
\lambda_{\frak s} = \ve_{1}+\ve_{p+1}& \ve_{1}+\ve_{p+2} &\cdots &\ve_{1}+\ve_{p+q-1}&\ve_{1}\pm \ve_{p+q}&\ve_{1}-\ve_{p+q-1}&\cdots&\ve_{1}-\ve_{p+1}\\
 \ve_{2}+\ve_{p+1}& \ve_{2}+\ve_{p+2} &\cdots &\ve_{2}+\ve_{p+q-1}&\ve_{2}\pm\ve_{p+q}&\ve_{2}-\ve_{p+q-1}&\cdots&\ve_{2}-\ve_{p+1}\\
 \vdots& \vdots& \vdots& \vdots& \vdots& \vdots& \vdots& \vdots\\
  \ve_{p-1}+\ve_{p+1}& \ve_{p-1}+\ve_{p+2} &\cdots &\ve_{p-1}+\ve_{p+q-1}&\ve_{p-1}\pm\ve_{p+q}&\ve_{p-1}-\ve_{p+q-1}&\cdots&\ve_{p-1}-\ve_{p+1}\\\hline
\pm \ve_{p}+\ve_{p+1}& \pm\ve_{p}+\ve_{p+2} &\cdots &\pm\ve_{p}+\ve_{p+q-1}& \pm\ve_{p}\pm\ve_{p+q}&\pm\ve_{p}-\ve_{p+q-1}&\cdots&\pm\ve_{p}-\ve_{p+1}\\\hline
 -\ve_{p-1}+\ve_{p+1}&  -\ve_{p-1}+\ve_{p+2} &\cdots & -\ve_{p-1}+\ve_{p+q-1}& -\ve_{p-1}\pm\ve_{p+q}& -\ve_{p-1}-\ve_{p+q-1}&\cdots& -\ve_{p-1}-\ve_{p+1}\\
 \vdots& \vdots& \vdots& \vdots& \vdots& \vdots& \vdots& \vdots\\
  -\ve_{1}+\ve_{p+1}& -\ve_{1}+\ve_{p+2} &\cdots &-\ve_{1}+\ve_{p+q-1}&-\ve_{1}\pm\ve_{p+q}&-\ve_{1}-\ve_{p+q-1}&\cdots&-\ve_{1}-\ve_{p+1}\\\hline
\end{array}
\end{equation}}
\caption{{Hasse diagram for} $\frak{so}(2p,2q)$}
\label{eveneven}
\end{figure}

The set itself is the same as in the case of $(m,n) = (2p+1, 2q+1)$, but the indices are different;
for ${a\leq p < b}$,
 \begin{eqnarray*}
 \ind(\ve_a -\ve_b)&=& (a-1)+(p+q-b) +q-1,  \\
 \ind(-\ve_a +\ve_b) &=& (b-(p+1)) + (p-a) +p -1.
 \end{eqnarray*}
Thus  $\ind(-\Lambda_{\ext}(E_0))=p-1$. \\

 Since we have $W_1^{\theta} =\{ w_{\bf j}, w'_{\bf j}:{\bf j} \subset \{1, \ldots, r\}\}$, where
  \begin{eqnarray*}
  w_{\bf j}^{-1}(\ve_i) &=&\ve_{j_i} \text{ for all } i\\
  {w_{\bf j}'}^{-1}(\ve_i) &=& \left\{
  \begin{array}{cc}
  \ve_{j_i} & \text{ for } i \not=p\\
  -\ve_{j_i} & \text{ for }i=p,r.
  \end{array} \right.
  \end{eqnarray*}
For ${\bf j} \subset \{1, \dots, r\}$, define $(a',b')$ and $(a'',b'')$ as in the case when $(m,n)=(2p+1, 2q+1)$. Then
  $\ind(-\Lambda_{\ext}(E_0^{w_{\bf j}}))=\min\{(b'-(p+1)) + (p-a') +p-1, (a'' -1) + (p+q-b'')+q-1\}$.

\medskip
 \noindent \textsf{\textbf{Case 4.}}  If $(m,n) = (4, 2q+1)$ with $q\geq 2$, {then $\frak g = \frak{so}(2r+1, \mathbb C)$ with $r=2+q$ and $\frak{k} = \frak{so}(4,\mathbb C)\times\frak{so}(2q+1,\mathbb C)$. In this case
$\frak h =\frak t = \{ \text{diag}(\ve_1, \ldots,
\ve_r,0, -\ve_{r},\ldots, -\ve_1): \ve_j\in \mathbb C \text{ for all } j\}$.
The simple root system of $\frak k$ is given by $\Psi_{\frak k}^{(1)} \cup \Psi_{\frak k}^{(2)}$ where
\begin{equation}\nonumber
\begin{aligned}
\Psi_{\frak k}^{(1)} &= \left\{ \ve_1-\ve_2
\right\}
=: \left\{ \beta^{(1)}
\right\}\\
\widetilde\Psi_{\frak k}^{(1)} &= \left\{ \ve_1+\ve_2
\right\}
=: \left\{ \tilde\beta^{(1)}
\right\}\quad\text{ and }  \\
\Psi_{\frak k}^{(2)} &= \{{\ve_3-\ve_4}, \ldots, \ve_{r-1}-\ve_r, \ve_r\}=: \{ \beta_1^{(2)},\ldots, \beta_{q}^{(2)}\}.
\end{aligned}
\end{equation}
}

Since $\lambda_{\frak s} = \xi^{(1)} + \tilde\xi^{(1)} + \xi^{(2)}_1$, 
{the Hasse diagram of $W.\lambda_{\frak s} $ is shown in Figure~\ref{so(4,2q+1)}. We construct it using a similar approach as in Figure~\ref{evenodd} and Figure~\ref{partial Hasse so(2p,2q+1)}.}

\begin{figure}[h]
{  \footnotesize
\begin{equation}\nonumber
\arraycolsep=1.4pt\def\arraystretch{1.5}\begin{array}{|cccccc|}\hline
 \xi^{(1)} + \tilde\xi^{(1)}  +\ve_3& \ldots& \xi^{(1)} + \tilde\xi^{(1)}  + \ve_{q+2}&\xi^{(1)} + \tilde\xi^{(1)}  -\ve_{q+2}& \ldots&\xi^{(1)} + \tilde\xi^{(1)} -\ve_3\\\hline
\xi^{(1)} - \tilde\xi^{(1)}  +\ve_3& \ldots&  \xi^{(1)} - \tilde\xi^{(1)} + \ve_{q+2}& \xi^{(1)} - \tilde\xi^{(1)}-\ve_{q+2}& \ldots& \xi^{(1)} - \tilde\xi^{(1)}-\ve_3\\
-\xi^{(1)} + \tilde\xi^{(1)} +\ve_3& \ldots&  -\xi^{(1)} + \tilde\xi^{(1)}  + \ve_{q+2}& -\xi^{(1)} + \tilde\xi^{(1)} -\ve_{q+2}& \ldots&-\xi^{(1)} + \tilde\xi^{(1)}-\ve_3\\\hline
-\xi^{(1)} - \tilde\xi^{(1)}+\ve_3& \ldots&-\xi^{(1)} - \tilde\xi^{(1)} + \ve_{q+2}& -\xi^{(1)} - \tilde\xi^{(1)}-\ve_{q+2}& \ldots&-\xi^{(1)} - \tilde\xi^{(1)}-\ve_3\\\hline
\end{array}
\end{equation}}
\caption{{Hasse diagram for} $\frak{so}(4,2q+1)$}
\label{so(4,2q+1)}
\end{figure}

Figure~\ref{so(4,2q+1)} is equal to
{  \footnotesize
\begin{equation}\nonumber
{\arraycolsep=1.4pt\def\arraystretch{1.5}\begin{array}{|cccccc|}\hline
{\color{blue}\ve_1 +\ve_3}& \ldots& {\color{blue}\ve_1  + \ve_{q+2}}&{\color{blue}\ve_1 -\ve_{q+2}}& \ldots& {\color{blue}\ve_1 -\ve_3}\\\hline
{\color{blue}\ve_2+\ve_3}& \ldots& {\color{blue}\ve_2  + \ve_{q+2}}& {\color{blue}\ve_2 -\ve_{q+2}}& \ldots&{\color{blue}\ve_2-\ve_3}\\
 -\ve_2  +\ve_3& \ldots& -\ve_2 + \ve_{q+2}& -\ve_2-\ve_{q+2}& \ldots& -\ve_2-\ve_3\\\hline
-\ve_1+\ve_3& \ldots&-\ve_1 + \ve_{q+2}&-\ve_1-\ve_{q+2}& \ldots&-\ve_1-\ve_3\\\hline
\end{array}}
\end{equation}}
and the roots highlighted in blue color are positive.
We have $w_{\text{max}} \lambda_{\frak s} =  \ve_2  - \ve_3$.
Since $-w_{\text{max}}\lambda_{\frak s} = -\ve_2  + \ve_3$, we obtain that
$\text{ind}(-w_{\text{max}}\lambda_{\frak s})=1$.\\

From now on for simplicity we will denote, for example, the first row of Figure~\ref{so(4,2q+1)}
 by $\ve_1 + (3,\cdots, q+2) \quad \ve_1 - (q+2,\cdots, 3)$.\\

$\bullet$
Select $j,k$ with $1\leq j<k \leq q+1$. Then by the Weyl group action related to
$w_{jk}\colon(1,2,\ldots, q+2)\mapsto (j,k,1,2, \ldots, \hat j, \ldots, \hat k,\ldots, q+2)$,
Figure~\ref{so(4,2q+1)} changes to Figure~\ref{after}.
\begin{figure}[h]
 {  \scriptsize\begin{equation}\nonumber
{ \arraycolsep=1.4pt\def\arraystretch{2}
\begin{array}{|cc|}\hline
 \ve_j + {\color{blue}(1,2, \ldots,j-1, \hat j, j+1, \ldots,k-1, \hat k,k+1,\ldots, q+2)}&
\ve_j - {({\color{blue}q+2,\ldots, k+1,\hat k, k-1, \ldots, j+1},\hat j, j-1,\ldots, 2,1)}\\\hline
\ve_k+ {\color{blue}(1,2, \ldots,j-1, \hat j,j+1, \ldots,k-1, \hat k,k+1,\ldots, q+2)}&
 \ve_k - {({\color{blue} q+2,\ldots, k+1},\hat k, k-1, \ldots, j+1,\hat j, j-1,\ldots, 2,1)}\\
-\ve_k + {({\color{blue}1,2, \ldots,j-1, \hat j,j+1, \ldots,k-1}, \hat k,k+1,\ldots, q+2)}&
-\ve_k - (q+2,\ldots, k+1,\hat k, k-1, \ldots, j+1,\hat j, j-1,\ldots, 2,1)\\\hline
-\ve_j + {({\color{blue}1,2, \ldots,j-1}, \hat j,j+1, \ldots,k-1, \hat k,k+1,\ldots, q+2)}&
-\ve_j- (q+2,\ldots, k+1,\hat k, k-1, \ldots, j+1,\hat j, j-1,\ldots, 2,1)\\\hline
\end{array}}
\end{equation}}
\caption{{Hasse diagram for} $\frak{so}(4,2q+1)$ with $w_{jk}$ }
\label{after}
\end{figure}

\noindent{The weights highlighted in blue color are positive in $\frak b^{\text{ref},n}$.}
Hence
\begin{equation}\label{+}
\begin{aligned}
w_{\text{max}} \lambda_{\frak s} &= \ve_j-\ve_{j+1},
\quad \text{ind}(-w_{\text{max}}\lambda_{\frak s})=j+1
\quad\text{ when } j+1\neq k\\
w_{\text{max}} \lambda_{\frak s} &= \ve_k-\ve_{k+1},
\quad \text{ind}(-w_{\text{max}}\lambda_{\frak s})=k-1=j
\quad\text{ when } j+1=k.\\
\end{aligned}
\end{equation}

$\bullet$
Select $j,k$ with $1<j<k \leq q+1$. Then by the Weyl group action according to
$(1,2,\ldots, q+2)\mapsto (j,-k,1,2, \ldots, \hat j, \ldots, \hat k,\ldots, q+2)$, the image of Figure~\ref{so(4,2q+1)} by
$w^{-1}$ is equal to Figure~\ref{after}. Hence
$\text{ind}(-w_{\text{max}}\lambda_{\frak s})$ is the same with \eqref{+}.

\medskip

 \noindent \textsf{\textbf{Case 5.}} If $(m,n) = (4, 2q)$ with $q\geq 2$, then it is similar to the case $(m,n) = (4, 2q+1)$.
By the same calculation we obtain that
 $\text{ind}(-w_{\text{max}}\lambda_{\frak s})$ is the same with \eqref{+} when the Borel subgroup is
 defined by the action related to
$(1,2,\ldots, q+2)\mapsto (j,\pm k,1,2, \ldots, \hat j, \ldots, \hat k,\ldots, {\pm(q+2)})$.

 \medskip

 \noindent \textsf{\textbf{Case 6.}} If $(m,n)=(2, 2q+1)$ with $q\geq 1$,
 {then $\lambda_+ = \ve_1+\ve_2$ because $G_0$ is of Hermitian type.}
 {The Hasse diagram of $W.\lambda_{+}$ is given by
\begin{center}
\begin{tikzpicture}
    \matrix (a) [matrix of math nodes, column sep=0.35cm ]{
    \ve_2+\ve_1
     &\ve_3+\ve_1
     & \cdots
     & \ve_{q+1}+\ve_1
     &-\ve_{q+1}+\ve_1
     & \cdots
     & -\ve_2+\ve_1\\
    };
    \draw[->] (a-1-1)--(a-1-2);
    \draw[->] (a-1-2)--(a-1-3);
    \draw[->] (a-1-3)--(a-1-4);
    \draw[->] (a-1-4)--(a-1-5);
    \draw[->] (a-1-5)--(a-1-6);
    \draw[->] (a-1-6)--(a-1-7);
\end{tikzpicture}
\end{center}
and the Hasse diagram of $W.\lambda_{-}$ is  given by
\begin{center}
\begin{tikzpicture}
    \matrix (a) [matrix of math nodes, column sep=0.35cm ]{
    \ve_2-\ve_1
     &\ve_3-\ve_1
     & \cdots
     & \ve_{q+1}-\ve_1
     &-\ve_{q+1}-\ve_1
     & \cdots
     & -\ve_2-\ve_1.\\
    };
    \draw[->] (a-1-1)--(a-1-2);
    \draw[->] (a-1-2)--(a-1-3);
    \draw[->] (a-1-3)--(a-1-4);
    \draw[->] (a-1-4)--(a-1-5);
    \draw[->] (a-1-5)--(a-1-6);
    \draw[->] (a-1-6)--(a-1-7);
\end{tikzpicture}
\end{center}}
\noindent We have $w_{\text{max}} \lambda_{+} = -\ve_2 + \ve_1$.
Since $-w_{\text{max}}\lambda_{+} = \ve_2 - \ve_1\in \Lambda^+$, we obtain that
$\text{ind}(-w_{\text{max}}\lambda_{+})=0$.\\

$\bullet$
Select an element $j$ from $\{1,\ldots, q+1\}$ and define $w_j\in W^\theta$
induced from $(1,\ldots, q+1)\mapsto (j, 1,\ldots,j-1, \hat j , j+1, \ldots, q+1)$.
Then $w_j^{-1}(W.\lambda_{+})$ is given by
{  \footnotesize
\begin{equation}\nonumber
\begin{aligned}
& {\color{blue} \ve_1+\ve_j,  \ldots, \ve_{j-1}+\ve_j, \ve_{j+1} + \ve_j, \ldots, \ve_{q+1} + \ve_j},  {\color{blue} -\ve_{q+1}+\ve_j, -\ve_q+\ve_j, \ldots, -\ve_{j+1} + \ve_j}, -\ve_{j-1}+\ve_j, \ldots, -\ve_1 + \ve_j
\end{aligned}
\end{equation}}
and $w_j^{-1}(W.\lambda_{-})$ is given by
{  \footnotesize
\begin{equation}\nonumber
\begin{aligned}
{\color{blue} \ve_1-\ve_j, \ldots, \ve_{j-1}-\ve_j}, \ve_{j+1} - \ve_j, \ldots, \ve_{q+1} - \ve_j,
 -\ve_{q+1}-\ve_j, -\ve_q-\ve_j, \ldots, -\ve_{j+1} - \ve_j, -\ve_{j-1}-\ve_j, \ldots, -\ve_1 - \ve_j
\end{aligned}
\end{equation}}
where the blue roots are positive in {$\frak b^{\text{ref},n}$}.
Since we have
 {$w_{\text{max,+}} \lambda_{+} = -\ve_{j+1} + \ve_j$} and hence {$-w_{\text{max},+}\lambda_{+} =  \ve_{j+1} - \ve_j$, we obtain
$\text{ind}(-w_{\text{max},+}\lambda_{+})=j-1$.}
\\

$\bullet$
Select an element $j$ from $\{1,\ldots, q+1\}$ and define $w_j\in W^\theta$
induced from $(1,\ldots, q+1)\mapsto (-j, 1,\ldots,j-1, \hat j , j+1, \ldots, q+1)$.
Then {$w_j^{-1}(W.\lambda_{+})$} is given by
{  \footnotesize
\begin{equation}\nonumber
\begin{aligned}
{\color{blue} \ve_1-\ve_j,  \ldots, \ve_{j-1}-\ve_j}, \ve_{j+1} - \ve_j, \ldots, \ve_{q+1} - \ve_j,
-\ve_{q+1}-\ve_j, -\ve_q-\ve_j, \ldots, -\ve_{j+1} - \ve_j, -\ve_{j-1}-\ve_j, \ldots, -\ve_1 - \ve_j
\end{aligned}
\end{equation}}
and $w_j^{-1}(W.\lambda_{-})$ is given by
{  \footnotesize
\begin{equation}\nonumber
\begin{aligned}
{\color{blue} \ve_1+\ve_j,  \ldots, \ve_{j-1}+\ve_j, \ve_{j+1} + \ve_j, \ldots, \ve_{q+1} + \ve_j},
{\color{blue} -\ve_{q+1}+\ve_j, -\ve_q+\ve_j, \ldots, -\ve_{j+1} + \ve_j}, -\ve_{j-1}+\ve_j, \ldots, -\ve_1 + \ve_j\\
\end{aligned}
\end{equation}}
where the blue roots are positive in $\frak b^{\text{ref},n}$.
Since {$w_{\text{max},-} \lambda_- = -\ve_{j+1} + \ve_j$ and hence $-w_{\text{max},-}\lambda_{-} =  \ve_{j+1} - \ve_j$, we obtain
$\text{ind}(-w_{\text{max},-}\lambda_{-})=j-1$.}

\medskip
 \noindent \textsf{\textbf{Case 7.}}  If $(m,n)=(2, 2q)$ with $q\geq 1$,  then {$\lambda_{+} = \ve_2 + \ve_1$ because $G_0$ is of Hermitian type. The Hasse diagram of $W.\lambda_+$ is given by
 \begin{center}
\begin{tikzpicture}
    \matrix (a) [matrix of math nodes, column sep=0.35cm ]{
    &&&&\ve_{q+1}+\ve_1&&&\\
     \ve_2+\ve_1
     &\ve_3+\ve_1
     &\cdots
     &\ve_q+\ve_1
     &~
     &-\ve_q+\ve_1
     &\cdots
     &-\ve_2+\ve_1\\
     &&&& -\ve_{q+1}+\ve_1&&&\\
    };
    \draw[->] (a-2-1)--(a-2-2);
    \draw[->] (a-2-2)--(a-2-3);
    \draw[->] (a-2-3)--(a-2-4);
    \draw[->] (a-2-4)--(a-1-5);
    \draw[->] (a-2-4)--(a-3-5);
    \draw[->] (a-1-5)--(a-2-6);
    \draw[->] (a-3-5)--(a-2-6);
    \draw[->] (a-2-6)--(a-2-7);
    \draw[->] (a-2-7)--(a-2-8);
\end{tikzpicture}
 \end{center}
and the Hasse diagram of $W.\lambda_-$ is given by
 \begin{center}
\begin{tikzpicture}
    \matrix (a) [matrix of math nodes, column sep=0.35cm ]{
    &&&&\ve_{q+1}-\ve_1&&&\\
     \ve_2-\ve_1
     &\ve_3-\ve_1
     &\cdots
     &\ve_q-\ve_1
     &~
     &-\ve_q-\ve_1
     &\cdots
     &-\ve_2-\ve_1.\\
     &&&& -\ve_{q+1}-\ve_1&&&\\
    };
    \draw[->] (a-2-1)--(a-2-2);
    \draw[->] (a-2-2)--(a-2-3);
    \draw[->] (a-2-3)--(a-2-4);
    \draw[->] (a-2-4)--(a-1-5);
    \draw[->] (a-2-4)--(a-3-5);
    \draw[->] (a-1-5)--(a-2-6);
    \draw[->] (a-3-5)--(a-2-6);
    \draw[->] (a-2-6)--(a-2-7);
    \draw[->] (a-2-7)--(a-2-8);
\end{tikzpicture}
 \end{center}}
We have {$w_{\text{max},+} \lambda_{+} = -\ve_2 + \ve_1$.
Since $-w_{\text{max},+}\lambda_{+} = \ve_2 - \ve_1\in \Lambda^+$, we obtain that
$\text{ind}(-w_{\text{max},+}\lambda_{+})=0$.}\\

$\bullet$
Select an element $j$ from $\{1,\ldots, q+1\}$ and define $w_j\in W^\theta$
induced from $(1,\ldots, q+1)\mapsto (j, 1,\ldots,j-1, \hat j , j+1, \ldots, q+1)$.
Then $w_j^{-1}(W.\lambda_{+})$ is given by
{\footnotesize
\begin{equation}\nonumber
{\color{blue} \ve_1+\ve_j, \ve_2+\ve_j, \ldots, \ve_{j-1}+\ve_j, \ve_{j+1} + \ve_j, \ldots,
\begin{array}{c}
\ve_{q+1} + \ve_j\\
 -\ve_{q+1}+\ve_j
 \end{array}},
 {\color{blue} -\ve_q+\ve_j, \ldots, -\ve_{j+1} + \ve_j}, -\ve_{j-1}+\ve_j, \ldots, -\ve_1 + \ve_j
\end{equation}}
and $w_j^{-1}(W.\lambda_{-})$ is given by
{\footnotesize
\begin{equation}\nonumber
 {\color{blue} \ve_1-\ve_j, \ve_2-\ve_j, \ldots, \ve_{j-1}-\ve_j}, \ve_{j+1} - \ve_j, \ldots,
\begin{array}{c}
\ve_{q+1} - \ve_j\\
  -\ve_{q+1}-\ve_j
  \end{array},
 -\ve_q-\ve_j, \ldots, -\ve_{j+1} - \ve_j, -\ve_{j-1}-\ve_j, \ldots, -\ve_1 - \ve_j
\end{equation}}
where the blue roots are positive in $\frak b^{\text{ref},n}$.
{Since $w_{\text{max},+} \lambda_+ = -\ve_{j+1} + \ve_j$ and $-w_{\text{max},+}\lambda_+ =  \ve_{j+1} - \ve_j$, we obtain that
$\text{ind}(-w_{\text{max},+}\lambda_{+})=j-1$\\

$\bullet$
Similarly one obtains $\text{ind}(-w_{\text{max},-}\lambda_{-})=j-1$.}

\subsection{$\frak{g}_0 = \frak{sp}(p,q)$}
{For the Lie algebra data for $\frak{sp}(p,q)$, we refer the reader to \cite[Appendix C.3]{K}.}
{The geometric description of $G/B$ for $\frak{sp}(p,q)$ is similar to that of $\frak{su}(p,q)$ and $\frak{so}(m,n)$.}
Since {$\frak g = \frak{sp}(r, \mathbb C)$ with $r=p+q$} and $\frak k = \frak {sp}(p,\mathbb C)\times \frak{sp}(q,\mathbb C)$, the Cartan subalgebras of them are
$\frak h =\frak t= \{ \text{diag}(\ve_1,\ldots, \ve_r, -\ve_r,\ldots, -\ve_1):\ve_j\in \mathbb C \text{ for all } j\}$ by considering $\frak g$ as a subset of $ M_{2r}(\mathbb C)$. The reference Borel subalgebras are given by
$\frak b^{\text{ref}} = \frak{sp}({r},\mathbb C)\cap B^-_{2r}$ and $\frak b_{\frak k}^{\text{ref}} = ( \frak{sp}(p,\mathbb C)\times \frak{sp}(q,\mathbb C))\cap B^-_{2r}$.
The simple root system of $\frak k$ is given by $\Psi_{\frak k}^{(1)} \cup \Psi_{\frak k}^{(2)}$ where
\begin{equation}\nonumber
\begin{aligned}
\Psi_{\frak k}^{(1)} &= \{ \ve_1-\ve_2, \ldots, \ve_{p-1}-\ve_p, 2\ve_p\}\quad\text{ and } \\
\Psi_{\frak k}^{(2)} &= \{ \ve_{p+1}-\ve_{p+2}, \ldots, \ve_{p+q-1}-\ve_{p+q}, 2\ve_{p+q}\}.
\end{aligned}
\end{equation}
Since $\lambda_{\frak s} =\xi_1^{(1)} + \xi_1^{(2)} = \ve_1 + \ve_{p+1}$, the orbit $W.\lambda_{\frak s}$ is given by Figure~\ref{oddodd} and $\text{ind}(-\Lambda_{\text{ext}}(E_0^{w_{\bf j}}))$ is the same to that of $\frak{so}(2p+1, 2q+1)$.

\subsection{$\frak g_0 = \frak{sl}(m; \mathbb H)$}
{The complex vector space $\mathbb C^{2m}$ carries a quaternionic vector space structure $\mathbb H^m$ defined by $\mathcal J\colon v\mapsto J\overline v$ where $J=\left(\begin{array}{cc}
0&I_m\\
-I_m&0
\end{array}\right)$.
The Lie group $SL(m,\mathbb H)$ is the set of elements in $SL(2m,\mathbb C)$ satisfying $\mathcal J X = X\mathcal J$.
We have 
$\frak g = \frak{sl}(2m)$.
A full flag $(0\subset V_1\subset \cdots \subset V_{2m-1}\subset \mathbb C^{2m})$ is said to be $\mathcal J$-generic if $\dim(V_i\cap \mathcal J(V_j))$ is minimal for each $j=1,\ldots, 2m-1$.
In \cite[Proposition 3.14]{HW02}, it is proved that there is just one open $SL(m, \mathbb H)$-orbit in $G/B$ and it is the set of all $\mathcal J$-generic flags.
}

{We have $\frak k_0 = \frak{sp}(m)$ and
$\frak k=\frak{ sp}(m, \mathbb C)$.
Note that $\Sigma(\frak g,\frak h) = \{\ve_i-\ve_j : 1\leq i\neq j\leq 2m\}$, $\Sigma(\frak k,\frak t) = \{  \pm (\ve_i\pm\ve_j), \pm2\ve_j : 1\leq i< j \leq m\}$.
The Cartan subalgebras of $\frak g$ and $\frak k$ are given by $\frak h = \{\text{diag}(\ve_1, \ldots , \ve_{2m}):\sum \ve_i = 0\}$ and
$\frak t = \{\text{diag}(\ve_1,\ldots, \ve_m,-\ve_m,\ldots, -\ve_1): \ve_j\in \mathbb C \text{ for all }j\}$, respectively.
Since $\lambda_\frak{s} = 2\ve_1$, the Hasse diagram of  $W. \lambda_\frak{s}$ is given by
\begin{center}
\begin{tikzpicture}
    \matrix (a) [matrix of math nodes, column sep=0.35cm ]{
    2\ve_1
     &2\ve_2
     & \cdots
     & 2\ve_{m}
     &-2\ve_m
     &-2\ve_{m-1}
     & \cdots
     & -2\ve_1\\
    };
    \draw[->] (a-1-1)--(a-1-2);
    \draw[->] (a-1-2)--(a-1-3);
    \draw[->] (a-1-3)--(a-1-4);
    \draw[->] (a-1-4)--(a-1-5);
    \draw[->] (a-1-5)--(a-1-6);
    \draw[->] (a-1-6)--(a-1-7);
    \draw[->] (a-1-7)--(a-1-8);
\end{tikzpicture}
\end{center}

Since $w_\text{max} \lambda_\frak{s} = 2\ve_m$, we have $\text{ind}(-w_\text{max} \lambda_\frak{s}) = m$.
Since the Cartan involution $\theta\colon SL(n,\mathbb C)\to SL(m,\mathbb C)$ is defined by $\theta(X) = (\overline X^t)^{-1}$, we obtain $|W^\theta|/|W_{\frak k}|=1$.}

%************************* diagram 22q+1, 22q
%\begin{figure}[!h]

 %       \centering

%        \includegraphics[angle=0, width=1.5\textwidth]{diagram 22q+1.png}

 %       \caption{so(2,2q+1)}

 %       \label{fig3}

%\end{figure}

%\begin{figure}[!h]

 %       \centering

  %      \includegraphics[angle=0, width=1.5\textwidth]{diagram 22q.png}

   %     \caption{so(2,2q)}

   %     \label{fig4}

%\end{figure}

%\subsubsection{Other cases}

 %   \subsubsection{$\frak{so}^*(2r)$} Hermitian

 %$r \geq 3$

 %\subsubsection{$\frak{e}_{6, T_1D_5}$} Hermitian

 %\subsubsection{$\frak{e}_{7, T_1E_6}$} Hermitian

 %\subsubsection{$\frak{sp}(p,q)$} equal rank

 %\subsubsection{$\frak{sl}(m, \mathbb H)$} non-equal rank

 %\subsubsection{Exceptional non-Hermitian cases}

%equal rank \\

%nonequal rank\\

 \section{Ampleness of base cycles on period domains}  \label{sect:general flag manifolds}

\subsection{From $G/B$ to $G/Q$}

Let $Q$ be a parabolic subgroup of $G$ containing $B^{\rm ref}$ and $\pi: G/B \rightarrow G/Q$ be the projection map. Then base cycles in $G/Q$ are of the form $\pi(C^{w_1})$ where $w_1 \in W_1^{\theta}$ and $C^{w_1}=K.[w_1B^{\rm ref}]$. For $w_1 \in W_1^{\theta}$ the central fiber $E_0^{w_1}$ at $\pi(w_1(B^{\rm ref})) \in G/Q$ is
$$E_0^{w_1} = \frak s/\frak s \cap w_1(\frak q),$$
which is the same as $\frak s \cap w_1(\frak q^n)$ if $\theta(\frak q)=\frak q$, i.e., the corresponding flag domain is measurable. The computation  of ampleness in this case is similar to the previous  computation except that $\frak b^{{\rm ref},n}$ is replaced by $\frak q^n$.
We will do the same computation for base cycles corresponding to period domains.

\subsection{Period domains}
{In this section, we will consider the ampleness of period domains. For general information on period domains, we refer the readers to \cite{CMP, V02}.}

A  {\it $\mathbb Z$-Hodge structure of weight} $n$ is a free $\mathbb Z$-module $V_{\mathbb Z}$ with
\begin{enumerate}
\item[(I)] a decomposition $V_{\mathbb C}= \oplus_{r+s=n}V^{r,s}$ such that $V^{s,r} =\overline{V^{r,s}}$; or equivalently,
\item[(II)] a filtration $0 \subset F^n \subset \dots \subset F^1 \subset F^0=V_{\mathbb C}$
 such that $V_{\mathbb C} =F^r \oplus \overline{F^{n-r+1}}$,

\end{enumerate}
where $V_{\mathbb C}=V_{\mathbb Z} \otimes_{\mathbb Z} \mathbb C$.
The equivalence between (I) and (II) is given as follows.
%\begin{enumerate}
%\item[(I) $\Rightarrow$ (II)]
\begin{equation}\nonumber
\begin{aligned}
V^{r,s} &\mapsto& F^r=\oplus_{\mu \geq r}V^{\mu, n-\mu} \\
%
%\item[(II) $\Rightarrow$ (I)]
F^r &\mapsto& V^{r,s} =F^r \cap \overline{F^{n-r}} .
\end{aligned}
\end{equation}
%\end{enumerate}

 %$\mathbb Q$, $\mathbb R$ \\

Set 
$$h^{r,s} :=\dim V^{r,s} \quad \text{ and } \quad f^r:=\dim F^r=  h^{n,0} + h^{n-1,1}+ \cdots   + h^{r,n-r} . $$

A {\it polarized Hodge structure  $(V , {\mathsf Q})$ of weight} $n$ is given by a Hodge structure of weight $n$ together with a nondegenerate form ${\mathsf Q}:V_{\mathbb Z} \otimes V_{\mathbb Z} \rightarrow \mathbb Z$ such that ${\mathsf Q}(v,w) = (-1)^n{\mathsf Q}(w,v)$ and
\begin{eqnarray*}
{\mathsf Q}(V^{r,s}, V^{r',s'})=0 && \text{ if } r' \not=s, \\
i^{r-s}{\mathsf Q}(V^{r,s}, \overline{V^{r,s}}) >0.
\end{eqnarray*}
The first condition is equivalent to
 \begin{eqnarray*}
{\mathsf Q}(F^r, F^{n-r+1})=0
 \end{eqnarray*}

For example, the second condition is expressed as follows. If $n=5$, then $V=V^{5,0} + V^{ 4,1} + V^{ 3,2} + V^{ 2,3} + V^{ 1,4} + V^{0,5}$ and $i{\mathsf Q}$ is a Hermitian form with
\begin{eqnarray*}
i {\mathsf Q}(V^{ 5,0}, \overline{V^{ 5,0}}) &>&0, \\
i{\mathsf Q}(V^{4,1}, \overline{V^{4,1}}) &<&0,\\
i{\mathsf Q}(V^{3,2}, \overline{V^{3,2}})&>&0.
\end{eqnarray*}
  If $n=6$, then $V_{\mathbb C}=V^{6,0} + V^{5,1} + V^{4,2} + V^{3,3} + V^{2,4} + V^{1,5}+V^{0,6}$ and ${\mathsf Q}$ is a Hermitian form with
\begin{eqnarray*}
  {\mathsf Q}(V^{6,0}, \overline{V^{ 6,0}}) &<&0 ,\\
 {\mathsf Q}(V^{5,1}, \overline{V^{ 5,1}}) &>&0,\\
 {\mathsf Q}(V^{4,2}, \overline{V^{ 4,2}})&<&0,\\
{\mathsf Q}(V^{3,3}, \overline{V^{3,3}}) &>& 0.
\end{eqnarray*}

Let $G_0=Aut(V_{\mathbb R}, {\mathsf Q})$ and $D$ be the set of polarized Hodge structures $(V, {\mathsf Q})$ with given Hodge numbers $\{h^{r,s} \}_{r+s=n}$, called a ({\it polarized}) {\it period domain}. 
The compact dual $\check D$ of $D$ is the set of { partial flags}  $F^{\bullet}=\{F^n \subset F^{n-1} \subset \dots \subset F^0=V_{\mathbb C}\}$ of $V_{\mathbb C}$ with $\dim F^r=f^r$ and
$${\mathsf Q}(F^r, F^{n-r+1})=0.$$
The last condition is equivalent to $F^r = (F^{n-r+1})^{\perp}$. Thus a partial flag $F^{\bullet}=\{F^n \subset F^{n-1} \subset \dots \subset F^0\}$ in $\check{D}$ is determined by $F^n, F^{n-1}, \dots, F^{k+1}$ where $n=2k+1$ or $n=2k$.

 \subsection{$ \frak g_0=\frak{sp}(m, \mathbb R)$}

 %Let $V$ be a complex vector space of dimension $2m$ with symplectic form $\omega$.
% Let $Fl_{\omega}(f^n, \dots, f^1)$ denote the partial flag manifold consisting of partial flags $0 \subset F^n \subset F^{n-1} \subset \dots \subset F^1 \subset F^0 =V$ such that $\dim F^p=f^p$ and $F^{n-p+1} ={F^p}^{\perp}$ where $p=1, \dots, n$.

Let $D$ be a period domain parameterizing polarized Hodge structures $(V,{\mathsf Q})$ with Hodge numbers $\{h^{r,s}\}_{r+s=n}$ with $n=2k+1$ for some $k\in \mathbb N$.
In this case $\dim V_{\mathbb C}=h^{n,0} + h^{n-1,1} + \dots + h^{1, n-1} + h^{0,n}={ 2f^{k+1}}$  is even. Put $ m:= {\dim V_{\mathbb C} \over 2} =f^{k+1} $. Then $G_0=Sp(m, \mathbb R)$ with symplectic form $\omega:=i{\mathsf Q}$ and the period domain $D$ consists of partial flags 
$${ 0 \subset F^n \subset \dots \subset  F^{k+1} =(F^{k+1})^{\perp} \subset \dots \subset (F^n)^{\perp} \subset V_{\mathbb C}}$$
in $\check D=Fl_{\omega}(f^n, \dots, f^{k+1}{;V_{\mathbb C} } )$ with $\dim F^r=f^r $
such that, letting $V^{r,s}=F^r \cap \overline{F^{n-r}}$,
\begin{equation}\label{Q1}
 (-1)^{r-k+1}\omega(V^{r,s}, \overline{V^{r,s}})>0, \quad 0 \leq r \leq k.
\end{equation}
This implies
$$D=Sp(m,\mathbb R)/U(h^{n,0})\times U(h^{ n-1,1}) \times \dots \times U(h^{ k+1,k}) $$
and { $\check D=Sp(m, \mathbb C)/Q$, where 
%$P_{f^n, f^{n-1}, \dots, f^{k+1}}$ 
$Q$ is the parabolic subgroup of $Sp(m, \mathbb C)$ marked at the subset $\{\psi_{f^n}, \dots, \psi_{f^{k+1}}\}$ of the simple root system $\Psi=\{\psi_1, \psi_2, \dots, \psi_m\}$ of $(\frak g, \frak h)$. In other words, the Lie algebra $\frak q^n$  of the unipotent part of $Q$ is the sum of root spaces of root $\alpha=\sum n_i\psi_i$ such that $\sum_{p=k+1}^n n_{f^p} >0$. } 

Depending on the parity of $k$, by \eqref{Q1}, $\omega$ is positive definite or negative definite as follows:

{
\begin{eqnarray*}
\begin{array}{c|c|c|c|c|c} 
k=\frac{n-1}{2} & V^{2k+1,0} & V^{2k ,1} & \quad \cdots    \qquad & V^{k+2,k -1} & V^{k+1,k} \\ 
\hline 
 \text{odd} & - &+ & &- &+ \\
  \text{even}  & + & - & & - & + \\
\end{array} 
\end{eqnarray*}
}

%\begin{equation}\nonumber
%\underbrace{V^{2k+1, 0}}_{\substack{- \text{ if } k \text{ is odd}\\+ \text{ if } k \text{ is even}}},
%\underbrace{V^{2k, 1}}_{\substack{+ \text{ if } k \text{ is odd}\\ - \text{ if } k \text{ is even}}}, \ldots,
%\underbrace{V^{k+1, k}}_{\substack{+ \text{ if } k \text{ is odd}\\+ \text{ if } k \text{ is even}}}
%\end{equation}

\noindent 
Therefore, the period domain $D$ is the flag domain in flag manifold $\check D$ corresponding to the element $w_{\bf j} \in W_1^{\theta}$, where
\begin{eqnarray*}
{\bf j}&=&\{f^n+1, \dots, f^{n-1}, \quad \dots , \quad  f^{k+2}+1,\ldots, f^{k+1}\},\\
\{1, \dots, m\}\backslash {\bf j}  &=&\{1, \dots, f^n, \quad  \dots, \quad f^{k+3}+1, \ldots, f^{k+2}\} \quad \text{ if } k \text{ is odd},
\end{eqnarray*}
and
\begin{eqnarray*}
{\bf j}&=&\{1, \dots, f^n,  \quad \dots , \quad  f^{k+2}+1,\ldots, f^{k+1}\} ,\\
\{1, \dots, m\}\backslash {\bf j}  &=&\{f^n+1, \dots, f^{n-1}, \quad \dots,\quad f^{k+3}+1, \ldots, f^{k+2}\}\quad \text{ if } k \text{ is even}.
\end{eqnarray*}
Note that $2 \ve_i=2\psi_i + 2\psi_{i+1} + \dots + 2 \psi_{m-1} + \psi_m$ for $i=1, \dots, m$ and $\frak q^n$ is the sum of root spaces of root $\alpha=\sum n_i\psi_i$ such that $\sum_{p=k+1}^n n_{f^p} >0$.
 % Put $c:=f^{k+1} - f^{k}=h^{k,k+1}$.
{Thus $w_{\bf j}^{-1}(W.\lambda_{\frak s}) \cap \Lambda(\frak q^{n}) = w_{\bf j}^{-1}(W.\lambda_{\frak s}) \cap \Lambda(\frak b^{\text{ref},n}) $. By the same computation    as in Section \ref{sp(r)}
 we get}
\begin{equation}\nonumber
\begin{aligned}
\ind(-\Lambda_{\ext}(E^{w_{\bf j}}_0))
&=\left\{
\begin{array}{ccc}
\min\{h_e, m-h_{e} \} & \text{ with }  h_{e} = {\displaystyle \sum_{r\, even, r\geq k+1}h^{r,s}}& \text{ if } k \text{ is odd} \\
\min\{h_o, m-h_{o} \} & \text{ with }  h_{o} = {\displaystyle\sum_{r\, odd, r\geq k+1}h^{r,s}} & \text{ if } k \text{ is even}
\end{array}
\right.\\
&={\min (h_o, h_e).}
\end{aligned}
\end{equation}
{Remark that the index value here coincides with the lower bound of Corollary 4 of \cite{HHL} which is a certain invariant on the pseudoconcavity of $Sp(m,\mathbb R)$-flag domains in the Grassmannian of $m$-planes.}

 \subsection{$ \frak g_0= \frak{so}(2p, \ell) $}
Let $D$ be a period domain parameterizing polarized Hodge structures $(V_{\mathbb C},{\mathsf Q})$ with Hodge numbers $\{h^{r,s}\}_{r+s=n}$ with $n=2k$ for some $k\in \mathbb N$.
In this case $\dim V_{\mathbb C}=m_e+ m_o$, where $m_e=\sum_{r, even}h^{r,s}$, $m_o=\sum_{r, odd}h^{r,s}$.  
 Put $m:= \left[ \dim V_{\mathbb C} \over 2 \right] = f^{k+1} +\left[h^{k,k} \over 2 \right]$.
 Then $G_0=SO(m_e, m_o)$ with 
symmetric form ${\mathsf Q}$  
and  the period domain $D$ consists of  partial flags
$$0 \subset F^n \subset \dots \subset F^{k+1} \subset  F^k= (F^{k+1})^{\perp} \subset \dots \subset (F^n)^{\perp} \subset V_{\mathbb C}$$
of $V_{\mathbb C}$ such that, letting $V^{r,s} = F^r \cap \overline{F^{n-r}}$,  
$$
(-1)^{r-k} {\mathsf Q}(V^{r,s}, \overline V^{r,s})>0.
$$
%e.g., if $n=4$, then $(m_e, m_o)=(2 h^{0,4} + h^{2,2}, 2 h^{1,3})$, and, if $n=6$, then $(m_e, m_o) = (2 h^{0,6} + 2 h^{2,4}, 2 h^{1,5} + h^{3,3})$.
   %  such that $(F^n,\ldots, F^k)\in Fl_{{\mathsf Q}}(f^n,\ldots, f^k)$.
This implies
$$D=SO(m_e, m_o)/U(h^{n,0})\times U(h^{ n-1,1}) \times \dots \times U(h^{ k+1, k-1 })\times SO(h^{k,k}) $$
and $\check D=SO(m_e+m_o, \mathbb C)/ Q$, where   
$Q$ is the parabolic subgroup of $SO(m_e+ m_o, \mathbb C)$ marked at the subset $\{\psi_{f^n}, \dots, \psi_{f^{k+1}}\}$ of the simple root system $\Psi=\{\psi_1, \psi_2, \dots, \psi_m\}$ of $(\frak g, \frak h)$.  In other words, the Lie algebra $\frak q^n$  of the unipotent part of $Q$ is the sum of root spaces of root $\alpha=\sum n_i\psi_i$ such that $\sum_{p=k+1}^n n_{f^p} >0$. 

Depending on the parity of $k$, ${\mathsf Q}$ is positive or negative definite as follows: 

{
\begin{eqnarray*}
\begin{array}{c|c|c|c|c|c} 
k=\frac{n}{2} & V^{2k,0} & V^{2k-1,1} & \quad \dots    \qquad & V^{k+1,k-1} & V^{k,k} \\ 
\hline 
 \text{odd} & - &+ & &- &+ \\
  \text{even}  & + & - & & - & + \\
\end{array} 
\end{eqnarray*}
}

%\begin{equation}\nonumber
%\underbrace{ V^{2k,0}}_{\substack{- \text{ if } k \text{ is odd}\\ + \text{ if } k \text{ is even}}},
%\underbrace{ V^{2k-1,1}}_{\substack{+ \text{ if } k \text{ is odd}\\ - \text{ if } k \text{ is even}}}, \ldots,
%\underbrace{ V^{k+1,k-1}}_{\substack{-\text{ if } k \text{ is odd }\\ -\text{ if } k \text{ is even}}}, 
%\underbrace{ V^{k,k}}_{\substack{+\text{ if } k \text{ is odd }\\ + \text{ if } k \text{ is even}}}
%\end{equation}
%Remark that if $k$ is even, then $m_o$ is even
%and if $k$ is odd, then $m_e$ is even. Recall that $f^r=\dim F^r$ where $F^r=\oplus_{\mu \geq r}V^{\mu, n-\mu}$. 
%Let $m:= \left[ \dim V_{\mathbb C} \over 2 \right]$. 
\noindent 
Therefore, { the period domain $D$ is the flag domain in flag manifold $\check{D}$} corresponding to the element  $w_{\bf j} \in W_1^{\theta}$
%$w_1\in W^\theta_1$   given by 
where  
%if $k$ is odd,
\begin{eqnarray*}
{\bf j}&=&\{f^{n}+1, \dots, f^{n-1}, \quad \dots, \quad {f^{k+1}+1, \ldots, m}\},\\
\{1, \dots, m\}\backslash {\bf j}  &=&\{1, \dots, f^{n},  \quad \dots , \quad f^{k+2}+1,\ldots, f^{k+1}\} \quad \text{ if } k \text{ is odd},
\end{eqnarray*}
 and %if $k$ is even,
\begin{eqnarray*}
{\bf j}&=&\{1, \dots, f^{n},  \quad \dots, \quad   {f^{k+1}+1, \ldots, m}\},\\
\{1, \dots, m\}\backslash {\bf j}  &=&\{f^{n}+1, \dots, f^{n-1}, \quad \dots , \quad f^{k+2}+1,\ldots, f^{k+1}\} \quad \text{ if } k \text{ is even}.
\end{eqnarray*}
  %Recall that $f^r=\dim F^r$ where $F^r=\oplus_{\mu \geq r}V^{\mu, n-\mu}$.  

  Note that { an element  in  $w_{\bf j}^{-1}(W.\lambda_{\frak s})$ is} of the form $\pm\ve_{j_a} \pm \ve_{j_b}$, where {$j_a \in {\bf j}$ and $j_b \in \{1, \dots, m\} \backslash {\bf j}$ and $\pm$ independent}. {This  can be expressed as the sum $\sum n_i\psi_i$ of simple roots, and satisfies  $\sum_{p=k+1 }^n n_{f^p} >0$ if it is contained in $\Lambda(\frak b^{\rm ref,n})$.}  Since $\frak q^n$ is the sum of root spaces of root $\alpha=\sum n_i\psi_i$ such that $\sum_{p=k+1}^n n_{f^p} >0$,  
   we have  $w_{\bf j}^{-1}(W.\lambda_{\frak s}) \cap \Lambda(\frak q^{n}) = w_{\bf j}^{-1}(W.\lambda_{\frak s}) \cap \Lambda(\frak b^{\text{ref},n}) $. 
\medskip

 If %$k$ is even , then $m_o$ is even
%and if 
 $k$ is odd (even, respectively), then $m_e$ ($m_o$, respectively) is even.    
We have $\frak g_0 = \frak{so}(2p, \ell)$ with
{ 
\begin{equation}\nonumber
\begin{aligned}
p &= h^{2k,0}+h^{2k-2,2} + \cdots + h^{k+1,k-1}= \frac{m_e}{2}, \\
\ell&=2h^{2k-1,1} + 2h^{2k-3,3} +\cdots+2h^{k+2,k-2} + h^{k,k} = m_o \text{ if } k \text{ is odd}, 
\end{aligned}
\end{equation}
 and 
\begin{equation}\nonumber
\begin{aligned}
p&=h^{2k-1,1} + h^{2k-3,3} +\cdots+h^{k+1,k-1} = \frac{m_o}{2},\\
\ell&=2h^{0,2k}+2h^{2k-2,2} + \cdots + 2h^{k+2,k-2} +  h^{k,k}=m_e \text{ if } k \text{ is even}.
\end{aligned}
\end{equation}
}

\noindent 
Let $q = \left[ \ell\over 2\right]$.
{Then $\frak g_0 =\frak{so}(2p,2q)$ or $\frak{so}(2p,2q+1)$. 
%$w_1^{-1} (\frak s)\cap \frak{q}^n $ is the set $\{ \ve_i\pm \ve_j : 1\leq i\leq p, \, p+1\leq j\leq p+q\}$.
 For these $p$ and $q$, let $(a', b')$ and $(a'', b'')$ be defined as in \eqref{ab}. By the same computation    as  in {\bf Case 2} and {\bf Case 3} of Section \ref{so(m,n)}, we have
 \begin{eqnarray*}
    && \ind(-\Lambda_{\ext}(E_0^{w_{\bf j}})) \\
     &=&
     \left\{\begin{array}{ll}
     \min\{(b'-(p+1)) + (p-a') +p-1, (a''-1)+ (p+q-b'') +q \} & \text{ if } \ell \text{ is odd},\\[5 pt]
      \min\{(b'-(p+1)) + (p-a') +p-1, (a''-1)+ (p+q-b'') +q-1 \}& \text{ if } \ell \text{ is even}. 
       \end{array} 
       \right.
 \end{eqnarray*}
}

%Let $q = \left[ \ell \over 2 \right]$.
%Then $w_1^{-1}(\frak s)\cap \frak{q}^n $ is the set $\{ \ve_i\pm\ve_j : 1\leq i\leq p, \, p+1\leq j\leq p+q\}$ since
%$w_1^{-1} \frak q \cap \frak{so}(2p,\ell)$ is compact.
%Let $(a', b')$ and $(a'', b'')$ be as \eqref{ab}.
%Then for given $p$ and $q$, we have
% $$\ind(-\Lambda_{\ext}(E_0^{w_{\bf j}}))= \min\{I^+({\bf j})-1, I^-({\bf j})- \delta_{\ell, e}\}$$
%where $\delta_{\ell,e}=1$ if $\ell$ is even and %$\delta_{\ell,e}=0$ if $m$ is odd.\\

\bigskip
\noindent
{\bf Competing interests: }The authors declare that they have no conflict of interest.


\begin {thebibliography} {XXX}

\bibitem[BE89]{BE}  
Baston, Robert J.; Eastwood, Michael G. 
{\it The Penrose transform. Its interaction with representation theory.} 
Oxford Mathematical Monographs. Oxford Science Publications. The Clarendon Press, Oxford University Press, New York, 1989. xvi+213 pp.

\bibitem[CMP17]{CMP}
Carlson, James; M\"uller-Stach, Stefan; Peters, Chris 
{\it Period mappings and period domains}. 
Second edition of [MR2012297]. Cambridge Studies in Advanced Mathematics, 168. Cambridge University Press, Cambridge, 2017. xiv+562 pp. ISBN: 978-1-316-63956-6; 978-1-108-42262-8 

\bibitem[FHW06]{FHW}
G. Fels, A. Huckleberry, and J. Wolf
{\it Cycle spaces of flag domains. A complex geometric viewpoint. Progress in Mathematics}, {\bf 245}. Birkh\"auser Boston, Inc., Boston, MA, 2006.

%\bibitem[H22]{Ha}
%T. Hayama, \emph{Generic 1-connectivity of flag domains in Hermitian symmetric spaces}, J. Lie Theory 32 (2022), no. 2, 553--561.

\bibitem[HHL19]{HHL}
T. Hayama, A. Huckleberry and  Q. Latif, 
{\it Pseudoconcavity of flag domains: the method of supporting cycles}, Math. Ann. (2019) {\bf 375}:671--685

\bibitem[HHS18]{HHS} 
Hong, Jaehyun; Huckleberry, Alan; Seo, Aeryeong 
{\it Normal bundles of cycles in flag domains}. 
S\~{a}o Paulo J. Math. Sci. {\bf 12} (2018), no. 2, 278--289.


\bibitem[H13]{H13} Huckleberry, Alan 
{\it Hyperbolicity of cycle spaces and automorphism groups of flag domains}, Amer. J. Math. {\bf 135} (2013), 291--310.

\bibitem[HSB01]{HSB01}
Huckleberry, A. T.; Simon, A.; Barlet, D. 
{\it On cycle spaces of flag domains of $SL_n(\mathbb R)$.}
J. Reine Angew. Math. {\bf 541} (2001), 171--208. 

\bibitem[HW02]{HW02}
Huckleberry, Alan T.; Wolf, Joseph A. 
{\it Cycle spaces of real forms of $SL_n(\mathbb C)$.} 
Complex geometry (G\"ottingen, 2000), 111--133, Springer, Berlin, 2002. 

\bibitem[K02]{K}
Knapp, Anthony W. 
{\it Lie groups beyond an introduction}. 
Second edition. Progress in Mathematics, {\bf 140}. Birkh\"auser Boston, Inc., Boston, MA, 2002. xviii+812 pp. ISBN: 0-8176-4259-5 22-01


\bibitem[M13]{M} 
 Matsumura, Shin-ichi 
 {\it Asymptotic cohomology vanishing and a converse to the Andreotti-Grauert theorem on surfaces}. 
 Ann. Inst. Fourier (Grenoble) {\bf 63} (2013), no. 6, 2199–-2221.




 \bibitem[OV90]{OV}
Onishchik, A. L.; Vinberg, \'E. B. {\it Lie groups and algebraic groups}. Translated from the Russian and with a preface by D. A. Leites. Springer Series in Soviet Mathematics. Springer-Verlag, Berlin, 1990. xx+328 pp.


\bibitem[S86]   {Sn1}
D. Snow,
{\it On the ampleness of homogeneous vector bundles}, Trans. Amer. Math. Soc. {\bf 294}  (1986),  no. 2, 585--594.

\bibitem[S]    {Sn2}
D. Snow,
{\it Homogeneous vector bundles}, available at https://www3.nd.edu/~snow/

\bibitem[S78] {S1}
Sommese, Andrew John 
{\it Submanifolds of Abelian varieties}.
Math. Ann. {\bf 233} (1978), no. 3, 229–-256. 


\bibitem[S83] {S2} 
Sommese, Andrew John 
{\it A convexity theorem.}
Singularities, Part 2 (Arcata, Calif., 1981), 497–-505,
Proc. Sympos. Pure Math., {\bf 40}, Amer. Math. Soc., Providence, RI, 1983.



\bibitem[V02]{V02}
Voisin, Claire 
{\it Hodge theory and complex algebraic geometry I}. Translated from the French original by Leila Schneps. Cambridge Studies in Advanced Mathematics, {\bf 76}. Cambridge University Press, Cambridge, 2002. x+322 pp. ISBN: 0-521-80260-1 

\end {thebibliography}
\end{document}